\documentclass[reqno]{amsart}
\usepackage[T1]{fontenc}
\usepackage{amsfonts}
\usepackage{amsmath}
\usepackage{amsthm}
\usepackage{amssymb}
\usepackage{graphicx}
\usepackage{subcaption}
\usepackage{url}
\usepackage{geometry}
\usepackage{hyperref}
\newtheorem{theorem}[subsection]{Theorem}
\newtheorem{definition}{Definition}[section]
\newtheorem{lemma}{Lemma}[section]
\newtheorem{proposition}{Proposition}[section]

\newtheorem{remark}{Remark}[section]

\title[A general quasilinear elliptic problem]{A general quasilinear elliptic problem with variable exponents and Neumann boundary conditions for image processing}
\author{Maxim Bogdan}
\email{maxim.bogdan.n6h@student.ucv.ro}
\date{\today}

\begin{document}
		
	\maketitle

	\begin{center} \footnotesize{Department of Mathematics, University of Craiova, Al. I. Cuza Street, no. 13, 200585, Craiova, Romania}
	\end{center}

\begin{abstract}
	The aim of this paper is to state and prove existence and uniqueness results for a general elliptic problem with homogeneous Neumann boundary conditions, often associated with image processing tasks like denoising. The novelty is that we surpass the lack of coercivity of the Euler-Lagrange functional with an innovative technique that has at its core the idea of showing that the minimum of the energy functional over a subset of the space $W^{1,p(x)}(\Omega)$ coincides with the global minimum. The obtained existence result applies to multiple-phase elliptic problems under remarkably weak assumptions.
\end{abstract}

\bigskip


\smallskip


\section{Introduction}

 When we want to act over a noisy grayscale image represented by $u_0$ (a function on $\Omega$ with values between 0 (black) and 1 (white)) with an evolution process towards a final improved image we can use as a model a general parabolic problem of the following form (see \cite{Alaa1}, \cite{guo1} for more details):

\begin{equation}\tag{$P$}\label{eqpg}
	\begin{cases}\dfrac{\partial u}{\partial t}-\operatorname{div}\mathbf{a}(x,\nabla u)=f\big (x,u(t,x)\big ), & (t,x)\in (0,\infty)\times\Omega\\[3mm] \dfrac{\partial u}{\partial\nu}=0, & (t,x)\in (0,\infty)\times\partial\Omega\\[3mm] u(0,x)=u_0(x), & x\in\Omega\end{cases}
\end{equation}

 If we want to know about the images we can get in the limit case $t\to\infty$ then we need to solve for the steady-states of \eqref{eqpg}. This is the subject of the present paper. To be more precise, in this article we will study the following problem:

\begin{equation}\tag{$E$}\label{eqeg}
	\begin{cases}-\operatorname{div}\mathbf{a}\big (x,\nabla U(x)\big )=f\big (x,U(x)\big ), & x\in\Omega\\[3mm] \dfrac{\partial U}{\partial\nu}=0, & x\in \partial\Omega\\[3mm] 0\leq U(x)\leq 1, & x\in\Omega\end{cases}
\end{equation}

 The Neumann homogeneous problem was studied for the $p(x)$-laplacian in \cite{Fan4}, \cite{mih1}, \cite{mih2}, \cite{boureanu}, \cite{deng1}, \cite{kim2} or for more general operators in \cite{kim1} but in the left-hand side of the equation there was a term of the form $\lambda |U(x)|^{p(x)-2}U(x)$ which plays a crucial role in the coercivity of the energy functional associated with the problem. Much more research was done for the associated Dirichlet problem: see for example \cite{Kimm1}, \cite{Rad1} and the references therein. 

\bigskip

 The key contributions of this article are as follows: the introduction of hypothesis \textbf{(H7)} under which we establish the existence of solutions for problem \eqref{eqeg} (which, in particular, includes multiple-phase problems) as well as the existence and uniqueness of solution for \eqref{eqelambda}. To the best of my knowledge, the latter had not been previously solved due to the lack of coercivity in the energy functional. Additionally, we present some other new results as Proposition \ref{propmult}, Lemma \ref{lemlocglob}, Theorem \ref{3thmrot} and Lemma \ref{3lemgama}. It is also important to mention that all the results of the paper are proved in full details and with very good references, which makes the article well-suited for beginners in this field.
 
 \bigskip

 We start section $2$ by stating the hypotheses under which we analyze problem \eqref{eqeg} along with the necessary notations. We then prove in Proposition \ref{propmult} that these hypotheses hold for a general multiple-phase operator. After that we present in section $3$ some basic results with complete proofs that will be used later in the paper. In section $4$ we give a \textbf{weak minimum principle} and a \textbf{weak comparison principle} for elliptic problems with variable exponents, both of which play a crucial role in our analysis. Section 5 focuses on the new problem \eqref{eqelambda}, which arises when applying Rothe’s method to the parabolic problem \eqref{eqpg}. Section 6 begins with the study of an approximation problem of \eqref{eqelambda}, followed by the presentation of the paper’s most significant results. Here, we employ a novel variational technique to show that a minimizer within a subset of $W^{1,p(x)}(\Omega)$ is, in fact, a global minimizer and thus a critical point of the energy functional (see Lemma \ref{lemlocglob}). In Section 7, using the sub-supersolution method, we construct the minimal and maximal solutions for \eqref{eqeg}. Finally, we conclude the paper with potential directions for future research, followed by an appendix containing the auxiliary results applied in this study.

\section{Hypotheses and notations}

\noindent We work under the following hypotheses:

\begin{enumerate}
	\item[\textbf{(H1)}] $\Omega\subset\mathbb{R}^N,\ N\geq 2$ is an open, bounded and connected Lipschitz domain.
	
	\bigskip
	
	\item[$\bullet$] $\mathcal{U}:=\big\{U\in L^{\infty}(\Omega)\ |\ 0\leq U\leq 1\ \text{a.e. on}\ \Omega \big \}$.
	
	\item[$\bullet$] For any $\lambda>0$, $\mathcal{M}_{\lambda}(\Omega)=\big \{u\in L^{\infty}(\Omega)\ |\ 0\leq u\leq \lambda\ \text{a.e. on}\ \Omega\big \}$. In particular $\mathcal{U}=\mathcal{M}_{1}(\Omega)$.
	
	\bigskip
	
	\item[\textbf{(H2)}] $p:\overline{\Omega}\to (1,\infty)$ is a log-H\"{o}lder continuous variable exponent with the property that: $p^-\geq\dfrac{2N}{N+2}$ or $p:\overline{\Omega}\to (1,\infty)$ is a continuous exponent with $p^->\dfrac{2N}{N+2}$. In both situations we may write $W^{1,p(x)}(\Omega)\hookrightarrow L^2(\Omega)$ with the mention that in the second situation we have that this embedding is also compact.\footnote{The proof can be seen in \cite[Theorem 2.3]{Fan2}. Also take a look at \cite[Theorem 2.2(c)]{Dinca2} and Proposition 2.2 from \cite{Fan3}.}
	
	\bigskip
	
	\item[$\bullet$] Denote $p^{-}=\displaystyle\min_{x\in\overline{\Omega}}\ p(x) >1$ and $p^+=\displaystyle\max_{x\in\overline{\Omega}}\  p(x)<\infty$. Let also $p'(x)=\dfrac{p(x)}{p(x)-1}$ be the conjugate variable exponent and $p^*(x)=\begin{cases} \dfrac{Np(x)}{N-p(x)}, & p(x)<N\\[3mm] +\infty, & p(x)\geq N\end{cases}$ the Sobolev critical exponent of $p(x)$.
	
	\begin{remark} We have $+\infty\geq p^*(x)\geq 2$ and $1\leq p^*(x)'\leq 2$ for all $x\in\overline{\Omega}$.
	\end{remark}
	
	\bigskip
	
	\item[\textbf{(H3)}] $\Psi:\overline{\Omega}\times (0,\infty)\to (0,\infty)$ with $\Psi(\cdot, s)$ measurable for each $s\in (0,\infty)$ and $\Psi(x,\cdot)\in\operatorname{AC}_{\text{loc}}\big ((0,\infty)\big )$ for a.a. $x\in\Omega$.\footnote{Hypothesis \textbf{(H3)} is taken from \cite{Rad2} and \cite{Kimm1}.}

	\item[\textbf{(H4)}] $\lim\limits_{s\to 0^+} \Psi(x,s)s=0$ for each $x\in\Omega$.
	
	\item[\textbf{(H5)}] For a.a. $x\in\Omega$ we have that $(0,\infty)\ni s\longmapsto \Psi(x,s)s\in (0,\infty)$ is a strictly increasing function.\footnote{Remark that I do not use the usual Leray-Lions structural conditions.}
	
	\bigskip

	\item[$\bullet$] We take $\mathbf{a}:\overline{\Omega}\times\mathbb{R}^N\to\mathbb{R}^N,\ \mathbf{a}(x,\xi)=\Psi(x,|\xi|)\xi$ if $\xi\neq \mathbf{0}$ and $\mathbf{a}(x,\textbf{0})=\textbf{0}$. Remark that $\forall\ \xi\in\mathbb{R}^N$ we have that $\mathbf{a}(\cdot,\xi)$ is measurable and for a.a. $x\in\Omega$ we have that $\mathbf{a}(x,\cdot)$ is continuous. Thus $\mathbf{a}$ is a Carath\'{e}odory mapping.\footnote{For definition and properties of Carath\'{e}odory functions I recommend \cite[Chapter 4]{Alip}.}
	
	\bigskip
	
	\item[$\bullet$] Define $\Phi:\overline{\Omega}\times\mathbb{R}\to [0,\infty),\ \Phi(x,s)=\begin{cases} \Psi(x,|s|)|s|, & s\neq 0\\ 0, & s=0\end{cases}$. This is an even function in the second argument. From \textbf{(H4)} and \textbf{(H5)} we have that $\Phi(x,\cdot)$ is strictly increasing on $[0,\infty)$ for a.a. $x\in\Omega$ and moreover it is continuous at $s=0$. Also, from \textbf{(H3)}, we conclude that $\Phi(x,\cdot)$ is continuous on $\mathbb{R}$ for a.a. $x\in\Omega$ and $\Phi(\cdot,s)$ is measurable for all $s\in\mathbb{R}$ (being a product of measurable functions). Hence $\Phi$ is a Carath\'{e}odory function. Notice that: $|\mathbf{a}(x,\xi)|=\Phi(x,|\xi|)$ and $\mathbf{a}(x,\xi)\cdot\xi=\Phi(x,|\xi|)|\xi|,\ \forall\ x\in\overline{\Omega},\ \forall\ \xi\in\mathbb{R}^N$.
	
	\bigskip
	
	\item[\textbf{(H6)}] $L^{p(x)}(\Omega)\ni v\longmapsto \Phi\big(\cdot,v(\cdot)\big )\in L^{p'(x)}(\Omega)$. This is the same as saying that the Nemytsky operator of $\Phi$ is defined as follows $\mathcal{N}_{\Phi}:L^{p(x)}(\Omega)\to L^{p'(x)}(\Omega)$. An other equivalent version would be: $L^{p(x)}(\Omega)^N\ni\mathbf{v}\longmapsto\mathbf{a}\big (\cdot,\mathbf{v}(\cdot)\big )\in L^{p'(x)}(\Omega)^N$.
	
	\bigskip
	
	\item[$\bullet$] $A_0:\overline{\Omega}\times[0,\infty)\to [0,\infty),\ A_0(x,s)=\displaystyle\int_{0}^{s} \Phi(x,\tau)\ d\tau$. From the \textit{Fundamental theorem of calculus} we have for a.a. $x\in\Omega$ that $A_0(x,\cdot)\in C^1\big ([0,\infty)\big )$ and $\dfrac{\partial A_0}{\partial t}(x,s)=\Phi(x,s),\ \forall\ s\in [0,\infty)$. 
	
	\item[$\bullet$] $A:\overline{\Omega}\times\mathbb{R}^N\to [0,\infty),\ A(x,\xi)=A_0(x,|\xi|)=\displaystyle\int_{0}^{|\xi|}\Phi(x,s)\ ds$. Note that $A_0$ and $A$ are well-defined since $\Phi(x,\cdot)$ is continuous on $[0,\infty)$ for a.a. $x\in\Omega$.

	\item[$\bullet$] $\mathcal{A}:W^{1,p(x)}(\Omega)\to [0,\infty)$ given by $\mathcal{A}(u)=\displaystyle\int_{\Omega} A(x,\nabla u(x))\ dx$.
	
	\bigskip
	
	\item[\textbf{(H7)}] There exists some constants $\delta>0$ and $\tilde{\delta}\geq 0$ such that:
	
	\begin{equation}
		\mathcal{A}(v)=\int_{\Omega} A(x,\nabla v(x))\ dx\geq \delta\int_{\Omega} |\nabla v(x)|^{p(x)}\ dx-\tilde{\delta},\ \forall\ v\in W^{1,p(x)}(\Omega).
	\end{equation}
	
	\bigskip
	
	\item[\textbf{(H8)}] $f:\overline{\Omega}\times [0,1]\to\mathbb{R}$ is a measurable function.
	
	\item[\textbf{(H9)}] The function $\Omega\ni x\longmapsto f(x,0)$ is from $L^{\infty}(\Omega)$.
	
	\item[\textbf{(H10)}] $f(x,0)\geq 0$ and $f(x,1)\leq 0$ for a.a. $x\in\Omega$. 
	
	\item[\textbf{(H11)}] There are two constants $\lambda_0,\gamma>0$ such that for a.a. $x\in\Omega$, the function $[0,1]\ni s\mapsto f(x,s)+\lambda_0 s$ is strictly increasing and the function $[0,1]\ni s\mapsto f(x,s)$ is $\gamma$--Lipschitz.

\begin{remark} If we have set $\gamma>0$ then any $\lambda_0>\gamma$ does the job, because for any reals $1\geq s_1>s_2\geq 0$ we get: $f(x,s_1)+\lambda_0 s_1-f(x,s_2)-\lambda_0 s_2\geq -\gamma(s_1-s_2)+\lambda_0 (s_1-s_2)=(\lambda_0-\gamma)(s_1-s_2)>0$.
\end{remark}

\begin{remark} In particular from $|f(x,s)-f(x,0)|\leq \gamma |s|\leq \gamma$ for a.a. $x\in\Omega$ and a.a. $s\in [0,1]$ we obtain that $|f(x,s)|\leq \gamma+\Vert f(\cdot,0)\Vert_{L^{\infty}(\Omega)}$ for a.a. $x\in\Omega$ and a.a. $s\in [0,1]$. Thus, if $U\in\mathcal{U}$ then $f(\cdot, U(\cdot))\in L^{\infty}(\Omega)$ and $\Vert f(\cdot,U(\cdot))\Vert_{L^{\infty}(\Omega)}\leq \gamma+\Vert f(\cdot,0)\Vert_{L^{\infty}(\Omega)}$. 
\end{remark}

\begin{remark}
	From \textnormal{\textbf{(H8)}} and \textnormal{\textbf{(H11)}} it follows that $f:\overline{\Omega}\times [0,1]\to\mathbb{R}$ is a Carath\'{e}odory function.
\end{remark}

\end{enumerate}

\begin{definition} We introduce the following Lipschitz extension function $\overline{f}:\overline{\Omega}\times\mathbb{R}\to\mathbb{R}$ given by:
	
	\begin{equation}
		\overline{f}(x,s)=\sup_{\tau\in [0,1]}f(x,\tau)-\gamma |s-\tau|=\begin{cases} f(x,0)+\gamma s, & s\in (-\infty,0)\\ f(x,s), & s\in [0,1] \\ f(x,1)-\gamma(s-1), & s\in (1,\infty) \end{cases}.
	\end{equation}

\end{definition}

\begin{proposition}\label{propel} The following properties of $\overline{f}$ hold: 
	\begin{enumerate}
		\item[\textnormal{\textbf{(1)}}] $\overline{f}(x,\cdot):\mathbb{R}\to\mathbb{R}$ is $\gamma-$Lipschitz for a.a. $x\in\Omega$;
		
		\item[\textnormal{\textbf{(2)}}] $|\overline{f}(x,s)|\leq \gamma |s|+\Vert f(\cdot,0)\Vert_{L^{\infty}(\Omega)}$ for any $(x,s)\in\overline{\Omega}\times\mathbb{R}$;
		
		\item[\textnormal{\textbf{(3)}}] For any $\overline{\lambda}_0>\gamma$ the function $\mathbb{R}\ni s\longmapsto\overline{f}(x,s)+\overline{\lambda}s$ is strictly increasing for a.a. $x\in\Omega$;
		
		\item[\textnormal{\textbf{(4)}}] $\overline{f}:\overline{\Omega}\times\mathbb{R}\to\mathbb{R}$ is a Carath\'{e}odory function.
	\end{enumerate}
\end{proposition}

\begin{proof}\textbf{(1)} Consider some real numbers $s_1,s_2$. We have that: $\overline{f}(x,s_1)=\displaystyle\sup_{\tau\in [0,1]} f(x,\tau)-\gamma|s_1-\tau|=\displaystyle\sup_{\tau\in [0,1]} f(x,\tau)-\gamma|(s_2-\tau)-(s_2-s_1)|\leq\displaystyle\sup_{\tau\in [0,1]} f(x,\tau)-\gamma|s_2-\tau|+\gamma |s_2-s_1|=\overline{f}(x,s_2)+\gamma |s_2-s_1|$. So $\overline{f}(x,s_1)-\overline{f}(x,s_2)\leq \gamma |s_2-s_1|$. By swapping $s_1$ and $s_2$ we get that $\overline{f}(x,s_2)-\overline{f}(x,s_1)\leq \gamma |s_1-s_2|$. Thus, $|\overline{f}(x,s_2)-\overline{f}(x,s_1)|\leq\gamma |s_2-s_1|$, as desired.
	
\noindent\textbf{(2)} From \textbf{(1)} we have that: $|\overline{f}(x,s)|-\Vert f(\cdot,0)\Vert_{L^{\infty}(\Omega)}\leq |\overline{f}(x,s)|-|f(x,0)|\leq |\overline{f}(x,s)-f(x,0)|=|\overline{f}(x,s)-\overline{f}(x,0)|\leq \gamma |s|$. 
	
\noindent\textbf{(3)} Let $s_1>s_2$ be two real numbers. Then $\overline{f}(x,s_1)+\overline{\lambda}_0 s_1-\overline{f}(x,s_2)-\overline{\lambda}_0 s_2=\overline{f}(x,s_1)-\overline{f}(x,s_2)+\overline{\lambda}_0(s_1-s_2)\geq -\gamma (s_1-s_2)+\overline{\lambda}_0(s_1-s_2)=(\overline{\lambda}_0-\gamma)(s_1-s_2)>0$.

\noindent\textbf{(4)} From \textbf{(1)} we have that $\overline{f}(x,\cdot)$ is continuous (being Lipschitz). It is also easy to observe that for each $s\in\mathbb{R}$ the function $\overline{f}(\cdot,s)$ is measurable (because $f(\cdot,s)$ is measurable). Thus $\overline{f}$ is a Carath\'{e}odory function.
\end{proof}

\begin{enumerate}
	\item[$\bullet$] $F:\overline{\Omega}\times [0,1]\to\mathbb{R}$, $F(x,s)=\displaystyle\int_{0}^{s} f(x,\tau)\ d\tau$.
	
	\item[$\bullet$] $\mathcal{F}:W^{1,p(x)}(\Omega)\cap\mathcal{U}\to\mathbb{R},\ \mathcal{F}(u)=\displaystyle\int_{\Omega} F(x,u(x))\ dx$.
	
	\item[$\bullet$]  $\overline{F}:\overline{\Omega}\times\mathbb{R}\to\mathbb{R}$, $\overline{F}(x,s)=\displaystyle\int_{0}^{s} \overline{f}(x,\tau)\ d\tau=\begin{cases} f(x,0)+\frac{\gamma}{2}s^2& s<0 \\[3mm] F(x,s)=\displaystyle\int_{0}^s f(x,\tau)\ d\tau, &  s\in [0,1] \\[3mm] F(x,1)+f(x,1)(s-1)-\frac{\gamma}{2}(s-1)^2, & s>1 \end{cases}.$
	
	\item[$\bullet$] $\overline{\mathcal{F}}:W^{1,p(x)}(\Omega)\to\mathbb{R},\ \overline{\mathcal{F}}(u)=\displaystyle\int_{\Omega} \overline{F}(x,u(x))\ dx$. Observe that $\overline{\mathcal{F}}(u)=\mathcal{F}(u)$ for any $u\in W^{1,p(x)}(\Omega)\cap\mathcal{U}$.
\end{enumerate}

\section{Preliminary results}

\begin{proposition}\label{propmult} The above hypotheses are satisfied by:
	\begin{equation}\label{divop}\mathbf{a}(x,\nabla u)=w_1(x)|\nabla u|^{p_1(x)-2}\nabla u+w_2(x)|\nabla u|^{p_2(x)-2}\nabla u+\hdots+w_{\ell}(x)|\nabla u|^{p_{\ell}(x)-2}\nabla u
	\end{equation}
	
	\noindent where $p_1,p_2,\dots, p_{\ell}:\overline{\Omega}\to (1,\infty),\ \ell\geq 1$ are continuous variable exponents such that $p:\overline{\Omega}\to (1,\infty),\ p(x)=\max\{p_1(x),p_2(x),\hdots,p_{\ell}(x)\}$ satisfies the inequality $p^-:=\displaystyle\min_{x\in\overline{\Omega}} p(x)>\dfrac{2N}{N+2}$,\footnote{Or $p_1,p_2,\hdots,p_{\ell}$ are log-H\"{o}lder continuous exponents and $p^-\geq\dfrac{2N}{N+2}$.} and $w_1,w_2,\hdots,w_{\ell}\in L^{\infty}(\Omega)$ with the property that there is some constant $\omega>0$ such that $\displaystyle\min_{j\in\overline{1,\ell}}\underset{x\in\Omega}{\operatorname{ess\ inf}}\ w_j(x)\geq\omega$.
\end{proposition}

\begin{proof} We have in this particular case that:
	
	\begin{enumerate}
		\item[$\bullet$] $\Psi:\overline{\Omega}\times (0,\infty)\to (0,\infty),\ \Psi(x,s)=w_1(x)s^{p_1(x)-2}+w_2(x)s^{p_2(x)-2}+\dots+w_{\ell}(x)s^{p_\ell(x)-2}$. So $\lim\limits_{s\to 0^+} \Psi(x,s)s=0$ for any $x\in\overline{\Omega}$.
		
		\item[$\bullet$] $\Phi:\overline{\Omega}\times\mathbb{R}\to [0,\infty),\ \Phi(x,s)=w_1(x)|s|^{p_1(x)-1}+w_2(x)|s|^{p_2(x)-1}+\dots+w_{\ell}(x)|s|^{p_{\ell}(x)-1}$ is strictly increasing on $[0,\infty)$ for any $x\in\Omega$.
		
		\item[$\bullet$] Hypothesis \textbf{(H6)} holds, Indeed, taking any $v\in L^{p(x)}(\Omega)$, for every $j\in\overline{1,\ell}$ we also have that $v\in L^{p_{j}(x)}(\Omega)$, since $p(x)\geq p_j(x)$ for all $x\in\overline{\Omega}$. So $|v|^{p_j(x)-1}\in L^{\frac{p_j(x)}{p_j(x)-1}}(\Omega)=L^{p'_j(x)}(\Omega)\subset L^{p'(x)}(\Omega)$,  because $p'(x)\leq p'_j(x)$ for all $x\in\overline{\Omega}$. Summing up we get that $\Phi\big (\cdot,v(\cdot)\big )\in L^{p'(x)}(\Omega)$.
		
		\item[$\bullet$] $A:\overline{\Omega}\times\mathbb{R}^N\to [0,\infty)$ and  $A(x,\xi)=w_1(x)\dfrac{|\xi|^{p_1(x)}}{p_1(x)}+w_2(x)\dfrac{|\xi|^{p_2(x)}}{p_2(x)}+\dots+w_{\ell}(x)\dfrac{|\xi|^{p_\ell(x)}}{p_\ell(x)}$.
		
		\item[$\bullet$] $\mathcal{A}:W^{1,p(x)}(\Omega)\to\mathbb{R},\ \mathcal{A}(u)=\displaystyle\int_{\Omega}w_1(x)\dfrac{|\nabla u|^{p_1(x)}}{p_1(x)}+w_2(x)\dfrac{|\nabla u|^{p_2(x)}}{p_2(x)}+\dots+w_{\ell}(x)\dfrac{|\nabla u|^{p_\ell(x)}}{p_\ell(x)}\ dx$

	\noindent To test hypothesis \textbf{(H7)} let's denote for each $j\in\overline{1,\ell}$: $\Omega_j=\{x\in\overline{\Omega}\ |\ p(x)=p_j(x)\}\subset\overline{\Omega}$. It is obvious that all of them are measurable and they cover $\overline{\Omega}$ meaning that: $\Omega_1\cup\Omega_2\cup\dots\cup\Omega_{\ell}=\overline{\Omega}$.
		
		\begin{align*}
			\mathcal{A}(u)&=\int_{\Omega}w_1(x)\dfrac{|\nabla u|^{p_1(x)}}{p_1(x)}+w_2(x)\dfrac{|\nabla u|^{p_2(x)}}{p_2(x)}+\dots+w_{\ell}(x)\dfrac{|\nabla u|^{p_\ell(x)}}{p_\ell(x)}\ dx\\
			&\geq \omega\int_{\Omega}\dfrac{|\nabla u|^{p_1(x)}}{p_1(x)}+\dfrac{|\nabla u|^{p_2(x)}}{p_2(x)}+\dots+\dfrac{|\nabla u|^{p_\ell(x)}}{p_\ell(x)}\ dx\\
			&=\omega\left ( \int_{\Omega}\dfrac{|\nabla u|^{p_1(x)}}{p_1(x)}\ dx+\int_{\Omega}\dfrac{|\nabla u|^{p_2(x)}}{p_2(x)}\ dx+\dots+\int_{\Omega}\dfrac{|\nabla u|^{p_\ell(x)}}{p_\ell(x)}\ dx\right ) \\
			&\geq \omega\left (\int_{\Omega_1}\dfrac{|\nabla u|^{p_1(x)}}{p_1(x)}\ dx+\int_{\Omega_2}\dfrac{|\nabla u|^{p_2(x)}}{p_2(x)}\ dx+\dots+\int_{\Omega_{\ell}}\dfrac{|\nabla u|^{p_\ell(x)}}{p_\ell(x)}\ dx\right ) \\
			&=\omega\left (\int_{\Omega_1}\dfrac{|\nabla u|^{p(x)}}{p(x)}\ dx+\int_{\Omega_2}\dfrac{|\nabla u|^{p(x)}}{p(x)}\ dx+\dots+\int_{\Omega_{\ell}}\dfrac{|\nabla u|^{p(x)}}{p(x)}\ dx\right ) \\
			&\geq \omega\int_{\Omega_1\cup\Omega_2\dots\cup\Omega_{\ell}}\dfrac{|\nabla u|^{p(x)}}{p(x)}\ dx\\
			&\geq \omega\int_{\Omega}\dfrac{|\nabla u|^{p(x)}}{p(x)}\ dx\\
			&\geq \dfrac{\omega}{p^+}\int_{\Omega} |\nabla u|^{p(x)}\ dx,\ \delta:=\dfrac{\omega}{p^+},\ \tilde{\delta}:=0.
		\end{align*}
	\end{enumerate}
	
\end{proof}

\begin{remark}
	Hypotheses \textnormal{\textbf{(H8)}}, \textnormal{\textbf{(H9)}}, \textnormal{\textbf{(H10)}} and \textnormal{\textbf{(H11)}} are satisfied by the following general logistic source with variable exponents: $f:\overline{\Omega}\times [0,1]\to\mathbb{R},\ f(x,s)=\alpha s^{q_1(x)}\left [r(x)-p(x)s^{q_2(x)} \right ]$ where $\alpha>0$ is a constant, $q_1:\overline{\Omega}\to [1,\infty),\ q_2:\overline{\Omega}\to [0,\infty)$ are two measurable and bounded exponents, and $r,p\in L^{\infty}(\Omega)$ with $0\leq r\leq p$ a.e. on $\Omega$. 
\end{remark}

\begin{proposition}\label{4prop1} The following basic and well-known properties hold:
	\begin{enumerate}
		\item[\textnormal{\textbf{(1)}}] For any $x\in\overline{\Omega}$ and any $\xi_1,\xi_2\in\mathbb{R}^N$ we have that:
		
		\begin{equation}
			\big (\mathbf{a}(x,\xi_1)-\mathbf{a}(x,\xi_2)\big )\cdot (\xi_1-\xi_2)\geq 0,
		\end{equation}
		
		\noindent with equality iff $\xi_1=\xi_2$.
		
		\item[\textnormal{\textbf{(2)}}] For a.a. $x\in\Omega$ and for all $\xi\in\mathbb{R}^N$ we have that $\nabla_{\xi} A(x,\xi)=\mathbf{a}(x,\xi)$.


		\item[\textnormal{\textbf{(3)}}] For a.a. $x\in\Omega$ one has:
		
		\begin{equation}\label{4eqAa}
			\mathbf{a}(x,\xi_2)\cdot (\xi_2-\xi_1)\geq A(x,\xi_2)-A(x,\xi_1)\geq \mathbf{a}(x,\xi_1)\cdot (\xi_2-\xi_1),\ \forall\ \xi_1,\xi_2\in\mathbb{R}^N.
		\end{equation}
		
		\noindent Equality occurs only when $\xi_1=\xi_2$.
		
		\item[\textnormal{\textbf{(4)}}] For a.a. $x\in\Omega$ we have that $A(x,\cdot)$ is an even function that is also strictly convex. 
		
		\item[\textnormal{\textbf{(5)}}] $\mathcal{N}_{\Phi}$ is a continuous and bounded operator and moreover there is some non-negative function $a\in L^{p'(x)}(\Omega)$ and a constant $b\geq 0$ such that: 
		
		\begin{equation}
			\Phi(x,s)\leq a(x)+b|s|^{p(x)-1},\ \forall\ s\in [0,\infty),\ \text{for a.a.}\ x\in\Omega.
		\end{equation}
		\noindent Thus the following inequality is true:
		
		\begin{equation}
			|\mathbf{a}(x,v)|\leq a(x)+b|v|^{p(x)-1},\ \forall\ v\in\mathbb{R}^N,\ \text{for a.a.}\ x\in\Omega.
		\end{equation}

		\item[\textnormal{\textbf{(6)}}] $\mathcal{A}$ is well-defined, $\mathcal{A}\in C^1\big (W^{1,p(x)}(\Omega)\big )$ and:
		
		\begin{equation}
			\langle \mathcal{A}'(u), \phi\rangle=\int_{\Omega} \mathbf{a}(x,\nabla u(x))\cdot \nabla\phi(x)\ dx,\ \forall\ u,\phi\in W^{1,p(x)}(\Omega).
		\end{equation}
		
		\item[\textnormal{\textbf{(7)}}] $\overline{F}$ is a Carath\'{e}odory function. Moreover for a.a. $x\in\Omega$ we have that $\overline{F}(x,\cdot)\in C^1(\mathbb{R})$ and $\dfrac{\partial \overline{F}}{\partial s}(x,s)=\overline{f}(x,s)$ for a.a. $x\in\Omega$ and $\forall\ s\in\mathbb{R}$.
		
		\item[\textnormal{\textbf{(8)}}] $\overline{\mathcal{F}}$ is well-defined, $\overline{\mathcal{F}}\in C^1\big (W^{1,p(x)}(\Omega)\big )$ and:
		
		\begin{equation}
			\langle \overline{\mathcal{F}}'(u),\phi \rangle=\int_{\Omega} \overline{f}(x,u(x))\phi(x)\ dx,\ \forall\ u,\phi\in W^{1,p(x)}(\Omega).
		\end{equation}
	\end{enumerate}
\end{proposition}

\begin{proof} \textbf{(1)} Fix some $x\in\overline{\Omega}$ and take any $\xi_1,\xi_2\in\mathbb{R}^N$. At that point, we distinguish three cases:

\begin{enumerate}
	\item[\textbf{(i)}] $\xi_1=\xi_2=\mathbf{0}$. Here there is nothing to prove.
	
	\item[\textbf{(ii)}] $\xi_1\neq \mathbf{0}$ and $\xi_2=\mathbf{0}$. We need to show that $\mathbf{a}(x,\xi_1)\cdot\xi_1>0\ \Leftrightarrow\ \Psi(x,|\xi_1|)\xi_1\cdot\xi_1>0\ \Leftrightarrow\ \Psi(x,|\xi_1|)|\xi_1|^2>0$ which is clearly true from \textbf{(H3)}.
	
	\item[\textbf{(iii)}] $\xi_1,\xi_2\neq\mathbf{0}$. We have that:

	\begin{align*}
		\big (\mathbf{a}(x,\xi_1)-\mathbf{a}(x,\xi_2)\big )\cdot (\xi_1-\xi_2)&=\Psi(x,|\xi_1|)|\xi_1|^2+\Psi(x,|\xi_2|)|\xi_2|^2-\big [\Psi(x,|\xi_1|)+\Psi(x,|\xi_2|) \big ]\xi_1\cdot\xi_2\\
		\text{(Cauchy ineq.)}	&\geq \Psi(x,|\xi_1|)|\xi_1|^2+\Psi(x,|\xi_2|)|\xi_2|^2-\big [\Psi(x,|\xi_1|)+\Psi(x,|\xi_2|) \big ]\cdot |\xi_1|\cdot|\xi_2|\\
		&=\big (|\xi_1|-|\xi_2| \big )\cdot \big (\Psi(x,|\xi_1|)|\xi_1|-\Psi(x,|\xi_2|)|\xi_2| \big )\geq 0\\
		&=\big (|\xi_1|-|\xi_2| \big )\cdot \big (\Phi(x,|\xi_1|)-\Phi(x,|\xi_2|) \big )\geq 0,
	\end{align*}
	
	\noindent because $\Phi(x,\cdot)$ is a strictly increasing function.
	
	\noindent Equality takes place only when $|\xi_1|=|\xi_2|$ and $\xi_1\cdot\xi_2=|\xi_1|\cdot |\xi_2|$ (i.e. there is some $\lambda\geq 0$ with $\xi_1=\lambda\xi_2$ or $\xi_2=\lambda\xi_1$). So we have that $\xi_1=\xi_2$ as claimed.
	
\end{enumerate}

\noindent\textbf{(2)} Let $x\in\overline{\Omega}$ such that $\Phi(x,\cdot)$ is continuous on $[0,\infty)$. If $\xi\neq 0$ we get from the \textit{chain rule} and the \textit{fundamental theorem of calculus} that:

\begin{equation}
	\nabla_{\xi} A(x,\xi)=\Phi(x,|\xi|)\cdot \nabla_{\xi}|\xi|=\Psi(x,|\xi|)\cdot|\xi|\cdot\dfrac{\xi}{|\xi|}=\Psi(x,|\xi|)\xi=\mathbf{a}(x,\xi). 
\end{equation}

\noindent For $\xi=0$ we get for each $k\in\overline{1,N}$ that:

\begin{equation}
	\dfrac{\partial A}{\partial\xi_k}(x,\mathbf{0})=\lim\limits_{\xi_k\to 0} \dfrac{A(x,\xi_k\mathbf{e}_k)-A(x,\mathbf{0})}{\xi_k}=\lim\limits_{\xi_k\to 0}\dfrac{1}{\xi_k}\int_{0}^{|\xi_k|} \Phi(x,s)\ ds=0.
\end{equation}

\noindent We have used the \textit{squeezing principle} and the monotonicity of $\Phi(x,\cdot)$: $\left | \dfrac{1}{\xi_k}\displaystyle\int_{0}^{|\xi_k|} \Phi(x,s)\ ds \right |\leq \dfrac{1}{|\xi_k|}\cdot |\xi_k|\cdot |\Phi(x,|\xi_k|)|=|\Phi(x,|\xi_k|)\stackrel{\xi_k\to 0}{\to} \Phi(x,0)=0$. Thus $\nabla_{\xi} A(x,\mathbf{0})=\mathbf{0}=\mathbf{a}(x,\mathbf{0})$.

\bigskip

\noindent\textbf{(3)} We prove just the right part of the inequality. The left part follows from the right one by simply switching $\xi_1$ and $\xi_2$. Suppose that $|\xi_2|\geq |\xi_1|$. If $\xi_1=\mathbf{0}$ we have that $\mathbf{a}(x,\xi_1)=\mathbf{0}$ and the inequality became $\displaystyle\int_0^{|\xi_2|}\Phi(x,s)\ ds\geq 0$ which is true, and the equality takes place only when $\xi_2=\mathbf{0}$.

\noindent For $|\xi_2|\geq |\xi_1|>0$ our inequality is equivalent to: $\displaystyle\int_{|\xi_1|}^{|\xi_2|} \Phi(x,s)\ ds\geq \Psi(x,|\xi_1|)\xi_1\cdot\xi_2-\Psi(x,|\xi_1|)|\xi_1|^2$. This is true since from the monotonicity of $\Phi(x,\cdot)$ we derive that $\displaystyle\int_{|\xi_1|}^{|\xi_2|} \Phi(x,s)\ ds\geq \displaystyle\int_{|\xi_1|}^{|\xi_2|} \Phi(x,|\xi_1|)\ ds=\big (|\xi_2|-|\xi_1| \big )\Psi(x,|\xi_1|)|\xi_1|=\Psi(x,|\xi_1|)|\xi_1|\cdot|\xi_2|-\Psi(x,|\xi|)|\xi_1|^2\geq \Psi(x,|\xi_1|)\xi_1\cdot\xi_2-\Psi(x,|\xi_1|)|\xi_1|^2$ from Cauchy inequality. Since $\Phi(x,\cdot)$ is strictly increasing the first inequality becomes equality only when $|\xi_1|=|\xi_2|$. The second inequality becomes equality iff $\xi_1\cdot\xi_2=|\xi_1|\cdot |\xi_2|$, i.e. $\xi_1=\lambda\xi_2$ or $\xi_2=\lambda\xi_1$ for some $\lambda>0$. Since $|\xi_1|=|\xi_2|$ it follows that $\lambda=1$ which means that $\xi_1=\xi_2$.

\bigskip

\noindent For $|\xi_1|\geq |\xi_2|$ and $\xi_2=\mathbf{0}$ the inequality is as follows: $-A(x,\xi_1)\geq -\mathbf{a}(x,\xi_1)\cdot \xi_1\ \Longleftrightarrow\ \Phi(x,|\xi_1|)|\xi_1|=\mathbf{a}(x,\xi_1)\cdot \xi_1\geq \displaystyle\int_{0}^{|\xi_1|}\Phi(x,s)\ ds$ which is true from the strict monotony of $\Phi(x,\cdot)$. Indeed: $\displaystyle\int_{0}^{|\xi_1|}\Phi(x,s)\ ds\leq \displaystyle\int_{0}^{|\xi_1|}\Phi(x,|\xi_1|)\ ds=\Phi(x,|\xi_1|)|\xi_1|$. Equality is attained only if $\xi_1=\mathbf{0}=\xi_2$.

\bigskip

\noindent For $|\xi_1|\geq |\xi_2|>0$ we need to show that: $A(x,\xi_1)-A(x,\xi_2)\leq \mathbf{a}(x,\xi_1)\cdot (\xi_1-\xi_2)$. The same arguments apply here:

\begin{align*}
	A(x,\xi_1)-A(x,\xi_2)&=\int_{|\xi_2|}^{|\xi_1|}\Phi(x,s)\ ds\leq \int_{|\xi_2|}^{|\xi_1|}\Phi(x,|\xi_1|)\ ds=\Phi(x,|\xi_1|)|\xi_1|-\Phi(x,|\xi_1|)|\xi_2|\\
	\text{(Cauchy ineq.)}\ \ \ \ \ \	 &=\mathbf{a}(x,\xi_1)\cdot\xi_1-\Psi(x,|\xi_1|)|\xi_1|\cdot|\xi_2|\leq \mathbf{a}(x,\xi_1)\cdot\xi_1-\Psi(x,|\xi_1|)\xi_1\cdot\xi_2\\
	&=\mathbf{a}(x,\xi_1)\cdot\xi_1-\mathbf{a}(x,\xi_1)\cdot\xi_2=\mathbf{a}(x,\xi_1)\cdot (\xi_1-\xi_2).
\end{align*}

\noindent Equality takes place only when $|\xi_1|=|\xi_2|$ and $\xi_1\cdot\xi_2=|\xi_1|\cdot |\xi_2|$ which means that $\xi_1=\xi_2$.

\bigskip

\noindent\textbf{(4)} $A(x,-\xi)=A(x,\xi)$ for any $\xi\in\mathbb{R}^N$. To show that $A(x,\cdot)$ is strictly convex we take $\xi_1,\xi_2\in\mathbb{R}^N$, $\lambda\in (0,1)$ and then using \textbf{(3)} we obtain that:

\begin{equation}
	\begin{cases} \lambda A(x,\xi_1)-\lambda A(x,\lambda\xi_1+(1-\lambda)\xi_2)\geq \lambda(1-\lambda)\mathbf{a}(x,\lambda\xi_1+(1-\lambda)\xi_2)\cdot(\xi_1-\xi_2)\\[3mm] (1-\lambda) A(x,\xi_2)-(1-\lambda)A(x,\lambda\xi_1+(1-\lambda)\xi_2)\geq -\lambda(1-\lambda)\mathbf{a}(x,\lambda\xi_1+(1-\lambda)\xi_2)\cdot (\xi_1-\xi_2)\end{cases}
\end{equation}

\noindent Adding these two relations will result in $\lambda A(x,\xi_1)+(1-\lambda)A(x,\xi_2)\geq A(x,\lambda\xi_1+(1-\lambda)\xi_2)$. Equality takes place only when $\xi_1=\xi_2=\lambda\xi_1+(1-\lambda)\xi_2$.

\bigskip

\noindent\textbf{(5)} The proof is a direct application of Theorem \ref{athnem}.

\bigskip

\noindent\textbf{(6)} First we need to explain why is $\mathcal{A}(u)$ a real number for each $u\in W^{1,p(x)}(\Omega)$. Since $A(x,\nabla u(x))\geq 0$ for a.a. $x\in\Omega$ we conclude that the integral $\displaystyle\int_{\Omega}A(x,\nabla u(x))\ dx$ exists, but it can be $\infty$. Now, using the fact that $\Phi(x,\cdot)$ is increasing on $[0,\infty)$ for a.a. $x\in\Omega$ we get that:

\begin{equation}
	0\leq A(x,\nabla u(x))=\int_{0}^{|\nabla u(x)|}\Phi(x,s)\ ds\leq \int_{0}^{|\nabla u(x)|}\Phi(x,|\nabla u(x)|)\ ds=\Phi(x,|\nabla u(x)|)|\nabla u(x)|,\ \text{for a.a.}\ x\in\Omega
\end{equation}

\noindent Since $u\in W^{1,p(x)}(\Omega)$ we have that $|\nabla u|\in L^{p(x)}(\Omega)$, and from \textbf{(H6)} we deduce that $\Phi(\cdot,|\nabla u|)\in L^{p'(x)}(\Omega)$. Finally, using \textit{H\"{o}lder inequality} we deduce that $\Phi(\cdot,|\nabla u|)|\nabla u|\in L^{1}(\Omega)$ as needed. So:

\begin{equation}
	0\leq \mathcal{A}(u)\leq \int_{\Omega}\Phi(x,|\nabla u(x)|)|\nabla u(x)|\ dx<\infty.
\end{equation}

\bigskip

\noindent Secondly, since $\phi\in W^{1,p(x)}(\Omega)$ we get that $|\nabla\phi|\in L^{p(x)}(\Omega)$. Thus:

\begin{align*}\left | \int_{\Omega} \mathbf{a}(x,\nabla u)\cdot\nabla\phi\ dx\right |&\leq \int_{\Omega} |\mathbf{a}(x,|\nabla u|)|\cdot |\nabla\phi|\ dx = \int_{\Omega} |\Phi(x,|\nabla u|)|\cdot |\nabla\phi|\ dx \\
	\text{(H\"{o}lder)}\ \ \ \ \ \ \ 	&\leq 2\Vert \Phi(x,|\nabla u|)\Vert_{L^{p'(x)}(\Omega)}\cdot \Vert |\nabla\phi(x)|\Vert_{L^{p(x)}(\Omega)}<\infty.
\end{align*}

\noindent So the mapping $x\longmapsto a(x,\nabla u(x))\cdot\nabla\phi(x)$ is in $L^1(\Omega)$ for each $\phi\in W^{1,p(x)}(\Omega)$. Therefore we can define the linear operator $L_u:W^{1,p(x)}(\Omega)\to\mathbb{R}$ by:

\begin{equation}
	L_u(\phi)=\int_{\Omega} \mathbf{a}(x,\nabla u(x))\cdot\nabla\phi(x)\ dx.
\end{equation}

\noindent Next, we prove that $L_u\in W^{1,p(x)}(\Omega)^*$, i.e. it is a bounded linear operator. Indeed, for any $\phi\in W^{1,p(x)}(\Omega)$ we have that:

\begin{align*}
	|L_u(\phi)|&=\left |\int_{\Omega}\mathbf{a}(x,\nabla u(x))\cdot\nabla\phi(x)\ dx \right |\leq \int_{\Omega} |\mathbf{a}(x,\nabla u(x))|\cdot |\nabla\phi(x)|\ dx\\
	&=\int_{\Omega} \Phi(x,|\nabla u(x)|)\cdot |\nabla\phi(x)|\ dx\\
	\text{(H\"{o}lder ineq.)}\ \ \ \ \ \ \		&\leq 2\Vert \Phi(\cdot,|\nabla u|)\Vert_{L^{p'(x)}(\Omega)}\cdot \Vert |\nabla\phi|\Vert_{L^{p(x)}(\Omega)}\leq  2\Vert \Phi(\cdot,|\nabla u|)\Vert_{L^{p'(x)}(\Omega)}\cdot \Vert \phi\Vert_{W^{1,p(x)}(\Omega)}.
\end{align*}

\noindent Observe now that, for any sequence $(\phi_n)_{n\geq 1}\subset W^{1,p(x)}(\Omega)$ with $\lim\limits_{n\to\infty} \Vert \phi\Vert_{W^{1,p(x)}(\Omega)}=0$ we may write:

\begin{align*}
	|\mathcal{A}(u+\phi_n)-\mathcal{A}(u)-L_{u}(\phi_n)|&=\left |\int_{\Omega} A(x,\nabla u+\nabla\phi_n)-A(x,\nabla u)-\mathbf{a}(x,\nabla u)\cdot \nabla\phi_n\ dx \right |\\
	&\leq\int_{\Omega} |A(x,\nabla u+\nabla\phi_n)-A(x,\nabla u)-\mathbf{a}(x,\nabla u)\cdot \nabla\phi_n|\ dx\\
	\text{\eqref{4eqAa}}\ \ \ \ \ \ \ 	&=\int_{\Omega} A(x,\nabla u+\nabla\phi_n)-A(x,\nabla u)-\mathbf{a}(x,\nabla u)\cdot \nabla\phi_n\ dx\\
	\text{\eqref{4eqAa}}\ \ \ \ \ \ \ 	&\leq \int_{\Omega} \big (\mathbf{a}(x,\nabla u+\nabla\phi_n)-\mathbf{a}(x,\nabla u)\big )\cdot \nabla\phi_n\ dx\\
	\text{(Cauchy ineq. and \textbf{(H6)})}\ \ \ \ \ \ \ &\leq  \int_{\Omega} \underbrace{\big |\mathbf{a}(\cdot,\nabla u+\nabla\phi_n)-\mathbf{a}(\cdot,\nabla u)\big |}_{\in L^{p'(x)}(\Omega)}\cdot |\nabla\phi_n|\ dx\\
	\text{(H\"{o}lder ineq.)}\ \ \ \ \ \ \ &\leq 2\big \Vert |\mathbf{a}(\cdot,\nabla u+\nabla\phi_n)-\mathbf{a}(\cdot,\nabla u)|\big \Vert_{L^{p'(x)}(\Omega)}\cdot \big \Vert |\nabla \phi_n|\big \Vert_{L^{p(x)}(\Omega)}\\
	&\leq 2\big \Vert |\mathbf{a}(\cdot,\nabla u+\nabla\phi_n)-\mathbf{a}(\cdot,\nabla u)|\big \Vert_{L^{p'(x)}(\Omega)}\cdot \Vert \phi_n \Vert_{W^{1,p(x)}(\Omega)}.
\end{align*}

\noindent Consider now $v:=\nabla u\in L^{p(x)}(\Omega)^N$ and define $\Phi_v:\overline{\Omega}\times\mathbb{R}^N\to\mathbb{R},\ \Phi_{v}(x,h)=|\mathbf{a}(x,v+h)-\mathbf{a}(x,v)|,\ \forall\ x\in\overline{\Omega},\ \forall\ h\in\mathbb{R}^N$ which is clearly a Carath\'{e}odory function. The Nemytsky operator associated to $\Phi_v$ is $\mathcal{N}_{\Phi_v}:L^{p(x)}(\Omega)^N\to L^{p'(x)}(\Omega)$. Indeed, since for all $h\in L^{p(x)}(\Omega)^N$ we have that:

\begin{equation*}
	\mathcal{N}_{\Phi_{v}}(h)=|\mathbf{a}(\cdot,v+h)-\mathbf{a}(\cdot,v)|\leq |\mathbf{a}(\cdot,v+h)|+|\mathbf{a}(\cdot,v)|=\mathcal{N}_{\Phi}(|v+h|)+\mathcal{N}_{\Phi}(|v|)\in L^{p'(x)}(\Omega).
\end{equation*}
\noindent we deduce from Theorem \ref{athleb} \textbf{(10)} in the Appendix that $\mathcal{N}_{\Phi_{v}}(h)\in L^{p'(x)}(\Omega)$. So, using Remark \ref{arenem}, we finally get that $\mathcal{N}_{\Phi_v}$ is continuous and bounded. Now, for $h_n:=\nabla\phi_n\in L^{p(x)}(\Omega)^N$ we know that $\Vert h_n\Vert_{L^{p(x)}(\Omega)^N}=\Vert |h_n|\Vert_{L^{p(x)}(\Omega)}=\Vert |\nabla\phi_n|\Vert_{L^{p(x)}(\Omega)}\longrightarrow 0$. Thus

\begin{equation*}\big \Vert |\mathbf{a}(\cdot,\nabla u+\nabla\phi_n)-\mathbf{a}(\cdot,\nabla u)|\big \Vert_{L^{p'(x)}(\Omega)}=\Vert\mathcal{N}_{\Phi_v}(h_n)-\mathcal{N}_{\Phi_v}(0)\Vert_{L^{p'(x)}(\Omega)}\longrightarrow 0.
\end{equation*}

\noindent Consequently, putting it all together we obtain that:

\begin{equation*}
	0\leq\lim\limits_{n\to\infty} \dfrac{|\mathcal{A}(u+\phi_n)-\mathcal{A}(u)-L_{u}(\phi_n)|}{\Vert\phi_n\Vert_{W^{1,p(x)}(\Omega)}}\leq 2 \lim\limits_{n\to\infty} \big \Vert |\mathbf{a}(\cdot,\nabla u+\nabla\phi_n)-\mathbf{a}(\cdot,\nabla u)|\big \Vert_{L^{p'(x)}(\Omega)}=0.
\end{equation*}

\noindent This proves that $\mathcal{A}$ is \textit{Fr\'{e}chet differentiable} on $W^{1,p(x)}(\Omega)$ and $\langle\mathcal{A}'(u),\phi\rangle=L_u(\phi),\ \forall\ \phi\in W^{1,p(x)}(\Omega)$.

\begin{remark}\label{4rema} We have proved at \textnormal{\textbf{(6)}} that if $v_n\longrightarrow v$ in $L^{p(x)}(\Omega)^N$ then $\mathbf{a}(\cdot, v_n)\longrightarrow \mathbf{a}(\cdot,v)$ in $L^{p'(x)}(\Omega)^N$.
\end{remark}

\bigskip

\noindent It remains to show that the (nonlinear) operator $\mathcal{A}':W^{1,p(x)}(\Omega)\to W^{1,p(x)}(\Omega)^*,\ \mathcal{A}'(u):=L_u,\ \forall\ u\in W^{1,p(x)}(\Omega)$ is continuous. First consider any $u\in W^{1,p(x)}(\Omega)$. Then:

\begin{align*}
	\Vert \mathcal{A}'(u)\Vert_{W^{1,p(x)}(\Omega)^*}&=\sup_{\phi\in W^{1,p(x)}(\Omega)\setminus\{0\}}\dfrac{\big |\langle\mathcal{A}'(u),\phi\rangle\big |}{\Vert\phi\Vert_{W^{1,p(x)}(\Omega)}}\leq \sup_{\phi\in W^{1,p(x)}(\Omega)\setminus\{0\}}\dfrac{\displaystyle\int_{\Omega}\big |\mathbf{a}(x,\nabla u)\big |\cdot\big |\nabla\phi\big |\ dx}{\Vert\phi\Vert_{W^{1,p(x)}(\Omega)}}\\
	\text{(H\"{o}lder ineq.)}\ \ \ \ \ \ \ &\leq \sup_{\phi\in W^{1,p(x)}(\Omega)\setminus\{0\}}\dfrac{2\big \Vert|\mathbf{a}(\cdot,\nabla u)|\big\Vert_{L^{p'(x)}(\Omega)} \cdot\big \Vert|\nabla\phi|\big\Vert_{L^{p(x)}(\Omega)}}{\Vert\phi\Vert_{W^{1,p(x)}(\Omega)}}\leq 2\big \Vert|\mathbf{a}(\cdot,\nabla u)|\big\Vert_{L^{p'(x)}(\Omega)}.
\end{align*}

\noindent Let us now take a sequence $(u_n)_{n\geq 1}\subset W^{1,p(x)}(\Omega)$ with $u_n\longrightarrow u$ in $W^{1,p(x)}(\Omega)$. Observe that:

\begin{align*}
	\big |\langle\mathcal{A}'(u_n)-\mathcal{A}'(u),\phi\rangle \big |&=\left |\int_{\Omega} \big (\mathbf{a}(x,\nabla u_n)-\mathbf{a}(x,\nabla u) \big )\cdot \nabla \phi\ dx \right |\leq \int_{\Omega} \big |\mathbf{a}(x,\nabla u_n)-\mathbf{a}(x,\nabla u) \big |\cdot |\nabla\phi|\ dx\\
	\text{(H\"{o}lder ineq.)}\ \ \ \ \ \ \ &\leq 2\big \Vert |\mathbf{a}(\cdot,\nabla u_n)-\mathbf{a}(\cdot,\nabla u)|\big \Vert_{L^{p'(x)}(\Omega)}\cdot \big\Vert |\nabla\phi |\big \Vert_{L^{p(x)}(\Omega)}\\
	&\leq 2\big \Vert |\mathbf{a}(\cdot,\nabla u_n)-\mathbf{a}(\cdot,\nabla u)|\big \Vert_{L^{p'(x)}(\Omega)}\cdot \Vert \phi \Vert_{W^{1,p(x)}(\Omega)}.
\end{align*}

\noindent We have that $\nabla u_n\longrightarrow \nabla u$ in $L^{p(x)}(\Omega)^N$ and from Remark \ref{4rema} we get that:

\begin{align*}
	\big \Vert\mathcal{A}'(u_n)-\mathcal{A}'(u)\Vert_{W^{1,p(x)}(\Omega)^*}&=\sup_{\phi\in W^{1,p(x)}(\Omega)\setminus\{0\}}\dfrac{|\langle\mathcal{A}'(u_n)-\mathcal{A}'(u),\phi\rangle \big |}{\Vert \phi \Vert_{W^{1,p(x)}(\Omega)}}\leq 2\big \Vert |\mathbf{a}(\cdot,\nabla u_n)-\mathbf{a}(\cdot,\nabla u)|\big \Vert_{L^{p'(x)}(\Omega)}\stackrel{n\to\infty}{\longrightarrow} 0.
\end{align*}

\noindent Thus $\mathcal{A}\in C^1\big (W^{1,p(x)}(\Omega) \big)$.

\noindent\textbf{(7)} From \textit{Fubini's theorem}\footnote{See Theorem 5.2.2 from \cite{Cohn}.} we have that for any fixed $s\in\mathbb{R}$ the function $\Omega\ni x\mapsto\overline{F}(x,s)=\displaystyle\int_{0}^s\overline{f}(x,\tau)\ d\tau$ is measurable (in fact even in $L^1(\Omega)$). It follows directly from the \textit{Fundamental theorem of calculus} that $\overline{F}(x,\cdot)\in C^1(\mathbb{R})$ for a.a. $x\in\Omega$ and $\dfrac{\partial \overline{F}}{\partial s}(x,s)=\overline{f}(x,s)$, for a.a. $x\in\Omega$ and for all $s\in\mathbb{R}$. In conclusion $\overline{F}$ is a Carath\'{e}odory function.

\noindent\textbf{(8)}  Since $f$ is measurable we may write for each $u\in W^{1,p(x)}(\Omega)\hookrightarrow\ L^{2}(\Omega)$ that:

\begin{align*}
	|\overline{\mathcal{F}}(u)|&\leq \int_{\Omega}\int_{0}^{|u(x)|} |\overline{f}(x,\tau)|\ d\tau\ dx\\
	\text{(Proposition \ref{propel} \ \textbf{(2)})}\ \ \ \ \ \ \ 	&\leq\int_{\Omega}\int_0^{|u(x)|} \gamma \tau+\Vert f(\cdot,0)\Vert_{L^{\infty}(\Omega)}\ d\tau\ dx\\
	&=\dfrac{1}{2}\int_{\Omega} u^2(x)\ dx+\Vert f(\cdot,0)\Vert_{L^{\infty}(\Omega)}\int_{\Omega} |u(x)|\ dx\\
\text{(Cauchy ineq.)}\ \ \ \ \ \ \ 	&\leq \dfrac{1}{2}\Vert u\Vert_{L^2(\Omega)}^2+\Vert f(\cdot,0)\Vert_{L^{\infty}(\Omega)}\Vert u\Vert_{L^2(\Omega)}\sqrt{|\Omega|}<\infty.
\end{align*}

\noindent This shows that indeed $\overline{\mathcal{F}}:W^{1,p(x)}(\Omega)\to\mathbb{R}$ is well-defined, and its definition makes sense on $L^{2}(\Omega)$ in fact.

\noindent Next, we'll show that $\overline{\mathcal{F}}$ is G\^ateaux-differentiable on $W^{1,p(x)}(\Omega)$. Set a direction $\phi\in W^{1,p(x)}(\Omega)$ (so $\phi\in L^{2}(\Omega)$) and consider any sequence $(\varepsilon_n)_{n\geq 1}\subset\mathbb{R}^*$ that is convergent to $0$. Without loss of generality we can assume that $|\varepsilon_n|\leq 1,\ \forall\ n\geq 1$. 

\noindent Consider for each $n\geq 1$ the function: $g_n:\Omega\to\mathbb{R},\ g_n(x)=\dfrac{\overline{F}(x,u(x)+\varepsilon_n\phi(x))-\overline{F}(x,u(x))}{\varepsilon_n}$. From the classical 1-dimensional \textit{mean value theorem} we have that $g_n(x)=\dfrac{\partial \overline{F}}{\partial s}(x,u(x)+s_x\phi(x))\phi(x)=\overline{f}(x,u(x)+s_n(x)\phi(x))\phi(x)$ for some real number $s_n(x)$ with $|s_n(x)|\leq\varepsilon_n$. Using the continuity of $\overline{f}(x,\cdot)$ for a.a. $x\in\Omega$ we get that:

\begin{equation}
	\lim\limits_{n\to\infty} g_n(x)=\lim\limits_{n\to\infty}\overline{f}(x,u(x)+s_n(x)\phi(x))\phi(x)\stackrel{|s_n(x)|\leq \varepsilon_n}{=}\overline{f}(x,u(x))\phi(x),\ \text{for a.a.}\ x\in\Omega.
\end{equation}

\noindent This proves that $g_n$ has a pointwise limit. The next step will be to show that $|g_n|$ is bounded by a function from $L^1(\Omega)$ for any $n\geq 1$. Indeed:

\begin{align*}
	|g_n(x)|&=|\overline{f}(x,u(x)+s_n(x)\phi(x))|\cdot |\phi(x)|\leq \left ( \gamma |u(x)+s_n(x)\phi(x)|+\Vert f(\cdot,0)\Vert_{L^{\infty}(\Omega)}\right )\cdot |\phi(x)|\ dx\\	&\leq \gamma |u(x)|\cdot |\phi(x)|+\gamma|\varepsilon_n| \phi^2(x)+\Vert f(\cdot,0)\Vert_{L^{\infty}(\Omega)}|\phi(x)|\\
	&\leq \gamma |u(x)|\cdot |\phi(x)|+\gamma \phi^2(x)+\Vert f(\cdot,0)\Vert_{L^{\infty}(\Omega)}|\phi(x)|\in L^1(\Omega).
\end{align*}

\noindent From the \textit{Lebesgue dominated convergence theorem} we get that:

\begin{align*}
	\lim\limits_{n\to \infty}\dfrac{\overline{\mathcal{F}}(u+\varepsilon_n\phi)-\overline{\mathcal{F}}(u)}{\varepsilon_n}&=\lim\limits_{n\to\infty} \int_{\Omega}\dfrac{\overline{F}(x,u(x)+\varepsilon_n\phi(x))-\overline{F}(x,u(x))}{\varepsilon_n}\ dx\\
	&=\lim\limits_{n\to\infty}\int_{\Omega} g_n(x)\ dx=\int_{\Omega} \overline{f}(x,u(x))\phi(x)\ dx.
\end{align*}

\noindent It is easy to remark that $\partial\overline{\mathcal{F}}(u):W^{1,p(x)}(\Omega)\to\mathbb{R},\ \partial_{\phi}\overline{\mathcal{F}}(u)=\displaystyle\int_{\Omega} \overline{f}(x,u(x))\phi(x)\ dx$ is a linear operator. We show now that it is in fact a bounded linear operator. Indeed:

\begin{align*}
	\left |\partial_{\phi}\overline{\mathcal{F}}(u) \right |&=\left | \int_{\Omega} \overline{f}(x,u(x))\phi(x)\ dx \right |\leq \int_{\Omega} |\overline{f}(x,u(x))|\cdot |\phi(x)|\ dx\\
	&\leq \int_{\Omega} \left (\gamma |u(x)|+\Vert f(\cdot,0)\Vert_{L^{\infty}(\Omega)}\right )|\phi(x)|\ dx\\
\text{(Cauchy ineq.)}	\ \ \ \ \ &\leq \big\Vert \gamma |u(x)|+\Vert f(\cdot,0)\Vert_{L^{\infty}(\Omega)}\big \Vert_{L^2(\Omega)}\Vert \phi\Vert_{L^2(\Omega)}\\
\big (W^{1,p(x)}(\Omega)\hookrightarrow L^{2}(\Omega)\big )\ \ \ \ \ \ \ &\lesssim \big\Vert \gamma |u(x)|+\Vert f(\cdot,0)\Vert_{L^{\infty}(\Omega)}\big \Vert_{L^2(\Omega)}\Vert \phi\Vert_{W^{1,p(x)}(\Omega)}
\end{align*}

\noindent Henceforth: $\displaystyle\sup_{\phi\in W^{1,p(x)}(\Omega)\setminus\{0\}}\dfrac{\big |\partial_{\phi}\overline{\mathcal{F}}(u) \big |}{\Vert \phi\Vert_{W^{1,p(x)}(\Omega)}}<\infty$, which shows that $\partial\overline{\mathcal{F}}(u)\in\mathcal{L}(W^{1,p(x)}(\Omega);\mathbb{R})=W^{1,p(x)}(\Omega)^*$, for any $u\in W^{1,p(x)}(\Omega)$. This observation allows us to write that: $\partial\overline{\mathcal{F}}:W^{1,p(x)}(\Omega)\to W^{1,p(x)}(\Omega)^*$. The following step is to prove that $\partial\overline{\mathcal{F}}$ is a continuous mapping.

\noindent Consider a sequence $(u_n)_{n\geq 1}\subset W^{1,p(x)}(\Omega)$ and some $u\in W^{1,p(x)}(\Omega)$ such that $u_n\to u$ in $W^{1,p(x)}(\Omega)$. Since $W^{1,p(x)}(\Omega)\hookrightarrow L^{2}(\Omega)$ we also have that $u_n\to u$ in $L^{2}(\Omega)$.

\begin{align*}\big \Vert\partial\overline{\mathcal{F}}(u_n)-\partial\overline{\mathcal{F}}(u)\big\Vert_{W^{1,p(x)}(\Omega)^*}&=\sup_{\phi\in W^{1,p(x)}(\Omega)\setminus\{0\}}\dfrac{\big |\partial_{\phi}\overline{\mathcal{F}}(u_n)-\partial_{\phi}\overline{\mathcal{F}}(u)\big |}{\Vert\phi\Vert_{W^{1,p(x)}(\Omega)}}\\
	&=\sup_{\phi\in W^{1,p(x)}(\Omega)\setminus\{0\}}\dfrac{\left |\displaystyle\int_{\Omega} \big (\overline{f}(x,u_n(x))-\overline{f}(x,u(x)))\cdot\phi(x) dx\right |}{\Vert\phi\Vert_{W^{1,p(x)}(\Omega)}} \\
	&\leq \sup_{\phi\in W^{1,p(x)}(\Omega)\setminus\{0\}}\dfrac{\displaystyle\int_{\Omega} \big |\overline{f}(x,u_n(x))-\overline{f}(x,u(x))|\cdot|\phi(x)| dx}{\Vert\phi\Vert_{W^{1,p(x)}(\Omega)}}\\
	&\leq \sup_{\phi\in W^{1,p(x)}(\Omega)\setminus\{0\}} \dfrac{\displaystyle\int_{\Omega} \gamma |u_n(x)-u(x)|\cdot|\phi(x)| dx}{\Vert\phi\Vert_{W^{1,p(x)}(\Omega)}}\\
	\text{(Cauchy ineq.)}\ \ \ \ \ \ \ \	&\leq \sup_{\phi\in W^{1,p(x)}(\Omega)\setminus\{0\}}\dfrac{\Vert u_n-u\Vert_{L^{2}(\Omega)}\Vert \phi\Vert_{L^{2}(\Omega)}}{\Vert\phi\Vert_{W^{1,p(x)}(\Omega)}}\\
	\big (W^{1,p(x)}(\Omega)\hookrightarrow L^{2}(\Omega)\big )\ \ \ \ \ \ \ &\leq \sup_{\phi\in W^{1,p(x)}(\Omega)\setminus\{0\}}\dfrac{\Vert u_n-u\Vert_{L^{2}(\Omega)}\Vert \phi\Vert_{W^{1,p(x)}(\Omega)}}{\Vert\phi\Vert_{W^{1,p(x)}(\Omega)}}\\
	&=\Vert u_n-u\Vert_{L^{2}(\Omega)}\stackrel{n\to\infty}{\longrightarrow} 0.
\end{align*}

\noindent Finally, from Theorem \ref{athmgatfre}, we obtain that $\overline{\mathcal{F}}\in C^1\big (W^{1,p(x)}(\Omega)\big )$ and:

\begin{equation*}
	\langle \overline{\mathcal{F}}'(u),\phi \rangle=\displaystyle\int_{\Omega} \overline{f}(x,u(x))\phi(x)\ dx,\ \forall\ u,\phi\in W^{1,p(x)}(\Omega).
\end{equation*}

\end{proof}

\section{Weak minimum principle and comparison principle for elliptic problems in variable exponent spaces}

\begin{theorem}[\textbf{Weak minimum principle}]\label{weakmp} Let some open, bounded and connected Lipschitz domain $\Omega\subset\mathbb{R}^N$ and then consider $c\in L^{\infty}(\Omega)^+$. If $u\in W^{1,p(x)}(\Omega)$, where $p:\overline{\Omega}\to (1,\infty)$ is a continuous exponent with $p^{-}>\dfrac{2N}{N+2}$\footnote{We can also assume that $p^{-}\geq\dfrac{2N}{N+2}$ if $p$ is log-H\"{o}lder continuous.}, satisfies in the \textbf{weak sense} the following inequalities:
	
	\begin{equation}
		\begin{cases}-\operatorname{div}\mathbf{a}(x,\nabla u(x))+c(x)u(x)\geq 0, & x\in\Omega\\[3mm]
			\dfrac{\partial u}{\partial\nu}(x)\geq 0, & x\in\partial\Omega\end{cases}
	\end{equation}
	
	\noindent i.e.
	
	\begin{equation}\label{ch3comp}
		\int_{\Omega} \mathbf{a}(x,\nabla u(x))\cdot\nabla\phi(x)\ dx+\int_{\Omega} c(x)u(x)\phi(x)\ dx\geq 0,\ \forall\ \phi\in W^{1,p(x)}(\Omega)^+,
	\end{equation}
	
	\noindent then exactly one of the following two situations takes place:
	
	\begin{enumerate}
		\item[\textnormal{\bf{(1)}}] $u(x)\geq 0$ for a.a. $x\in\Omega$.
		
		\item[\textnormal{\bf{(2)}}] $u$ is a constant strict negative function and $c\equiv 0$.
	\end{enumerate}
\end{theorem}

\begin{proof} First let's mention that from $p(x)\geq p^{-}>\dfrac{2N}{N+2}$ a.e. on $\Omega$ we get that $p^*(x)>2$ a.e. on $\Omega$ and thenceforth $W^{1,p(x)}(\Omega)\hookrightarrow L^2(\Omega)$. This shows that $u,\phi\in W^{1,p(x)}(\Omega)\subset L^2(\Omega)\ \Rightarrow\ u\phi\in L^1(\Omega)$ (from \textit{Cauchy inequality}). So, taking into account that $c\in L^{\infty}(\Omega)$, we conclude that $cu\phi\in L^1(\Omega)$. This is important in order to have finite integrals in \eqref{ch3comp}.

	\noindent Let $u^+=\max\{u,0\}$ and $u^-=-\min\{u,0\}$. From Theorem \ref{apthplus} we have that $u^+,u^{-}\in W^{1,p(x)}(\Omega)^+$. Choosing $\phi=u^{-}\in W^{1,p(x)}(\Omega)^+$ in \eqref{ch3comp} we reach to:
	
	\begin{align*}
		0&\leq \int_{\Omega} \Psi(x,|\nabla u(x)|)\big (\nabla u^+-\nabla u^{-})\cdot\nabla u^{-}\ dx+\int_{\Omega} c(u^+-u^-)u^{-}\ dx\\
		&=-\int_{\Omega}\Psi(x,|\nabla u(x)|)|\nabla u^{-}|^2\ dx-\int_{\Omega} c(u^{-})^2\ dx=-\int_{\Omega}\Phi(x,|\nabla u(x)|)|\nabla u^{-}|\ dx-\int_{\Omega} c(u^{-})^2\ dx\leq 0.
	\end{align*}
	
	\noindent Thus: $\displaystyle\int_{\Omega}\Phi(x,|\nabla u(x)|)|\nabla u^{-}|\ dx+\int_{\Omega} c(u^{-})^2\ dx=0$, which means that $\Phi(x,|\nabla u(x)|)|\nabla u^{-}(x)|=0$ a.e. on $\Omega$. If $\Phi(x,|\nabla u(x)|)=0$ then $\nabla u(x)=\mathbf{0}$ and therefore $\nabla u^{-}(x)=\mathbf{0}$. In the other case $\nabla u^{-}(x)=\mathbf{0}$. So we may conclude that $\nabla u^{-}=\mathbf{0}$ a.e. on $\Omega$. Now, since $u^{-}\in W^{1,p(x)}(\Omega)\subset W^{1,p^-}(\Omega)$ and $\Omega$ is a connected domain, we deduce that $u^{-}$ is a constant positive function.\footnote{Here, it is essential that $\Omega$ is connected. Take a look at the following difficult to prove fact: \cite[Proposition 7.18]{Salsa}.} If $u^{-}\equiv 0$ we are done since $u=u^{+}\geq 0$ a.e. on $\Omega$. If $u^{-}(x)=k>0$ a.e. on $\Omega$ is a strict positive constant function then $u=-u^{-}=-k$ (a strict negative constant function). Moreover $0=\displaystyle\int_{\Omega} c(x)(u^{-}(x))^2\ dx=k^2\int_{\Omega} c(x)\ dx\ \Rightarrow\ \int_{\Omega} c(x)\ dx=0$ which together with $c\geq 0$ a.e. on $\Omega$ imply that $c\equiv 0$.
	
\end{proof}

\begin{theorem}[\textbf{Weak comparison principle}]\label{weakcp} Fix some open, bounded and connected Lipschitz domain $\Omega\subset\mathbb{R}^N$ and then take $c\in L^{\infty}(\Omega)^+$. Let now $u_1,u_2\in W^{1,p(x)}(\Omega)$, where $p:\overline{\Omega}\to (1,\infty)$ is a continuous exponent with $p^{-}>\dfrac{2N}{N+2}$, and $g_1,g_2\in L^{\infty}(\Omega)$ with $g_2(x)\geq g_1(x)$ for a.a. $x\in\Omega$. If we have in the \textbf{weak sense} that:
	
	\begin{equation}
		\begin{cases}-\operatorname{div}\mathbf{a}(x,\nabla u_2(x))+c(x)u_2(x)\geq g_2(x)\geq g_1(x)\geq -\operatorname{div}\mathbf{a}(x,\nabla u_1(x))+c(x)u_1(x), & x\in\Omega\\[3mm]
			\dfrac{\partial u_2}{\partial\nu}(x)\geq 0\geq \dfrac{\partial u_1}{\partial\nu}(x) , & x\in\partial\Omega \end{cases}
	\end{equation}
	
	\noindent meaning that:
	
	\begin{equation}\label{ch3eqint}
		\begin{cases}\displaystyle\int_{\Omega} \mathbf{a}(x,\nabla u_2(x))\cdot\nabla\phi(x)\ dx+\int_{\Omega} c(x)u_2(x)\phi(x)\ dx\geq\int_{\Omega} g_2(x)\phi(x)\ dx\\[5mm]\displaystyle\int_{\Omega} \mathbf{a}(x,\nabla u_1(x))\cdot\nabla\phi(x)\ dx+\int_{\Omega} c(x)u_1(x)\phi(x)\ dx\leq\int_{\Omega} g_1(x)\phi(x)\ dx \end{cases}\ \forall\ \phi\in W^{1,p(x)}(\Omega)^+,
	\end{equation}
	
	\noindent then exactly one of the following two situation will happen:
	
	\begin{enumerate}
		\item[\textnormal{\bf{(1)}}] $u_2(x)\geq u_1(x)$ for a.a. $x\in\Omega$.
		
		\item[\textnormal{\bf{(2)}}] $u_2-u_1$ is a constant strict negative function, $c\equiv 0$ and $g_1\equiv g_2$.
	\end{enumerate}
	
	\noindent Moreover, if $u_1\equiv u_2$ then $g_1\equiv g_2$.
\end{theorem}

\begin{proof} Multiplying the second inequality in \eqref{ch3eqint} by $-1$ and then adding the two inequalities will result in:
	
	\begin{equation}
		\int_{\Omega} \big (\mathbf{a}(x,\nabla u_2(x))-\mathbf{a}(x,\nabla u_1(x))\big )\cdot\nabla\phi(x)\ dx+\int_{\Omega} c(u_2-u_1)\phi\ dx\geq \int_{\Omega}(g_2-g_1)\phi\ dx,\ \forall\ \phi\in W^{1,p(x)}(\Omega)^+.
	\end{equation}
	
	\noindent Choose the test function $\phi=(u_2-u_1)^{-}\in W^{1,p(x)}(\Omega)^+$ and denote $\Omega_1=\big\{x\in\Omega\ |\ u_2(x)<u_1(x)\big\}$. Since $u_2,u_1$ are measurable functions we have that $\Omega_1\subset\Omega$ is a bounded measurable set. Hence we can write:
	
	\begin{equation}\label{3ineqcomp}
		\int_{\Omega}\big (\mathbf{a}(x,\nabla u_2(x))-\mathbf{a}(x,\nabla u_1(x))\big )\cdot\nabla (u_2-u_1)^-\ dx-\int_{\Omega} c[(u_2-u_1)^-]^2\ dx\geq \int_{\Omega} (g_2-g_1)(u_2-u_1)^{-}\ dx.
	\end{equation}
	
	\noindent Note that $\nabla (u_2-u_1)^{-}=\begin{cases}\mathbf{0},\ \text{a.e. on}\ \Omega\setminus\Omega_1\\[3mm]-(\nabla u_2-\nabla u_1),\ \text{a.e. on}\ \Omega_1 \end{cases}$. Going further we get that:
	
	\begin{equation}
		0\geq-\int_{\Omega_1} \big (\mathbf{a}(x,\nabla u_2(x))-\mathbf{a}(x,\nabla u_1(x))\big )\cdot\big (\nabla u_2-\nabla u_1\big )\ dx-\int_{\Omega} c[(u_2-u_1)^-]^2\ dx\geq \int_{\Omega} (g_2-g_1)(u_2-u_1)^{-}\ dx\geq 0.
	\end{equation}
	
	\noindent Here we have used the strict monotony of $\mathbf{a}$ i.e. Proposition \ref{4prop1} \textbf{(1)}. So we can infer that $\big (\mathbf{a}(x,\nabla u_2(x))-\mathbf{a}(x,\nabla u_1(x))\big )\cdot\big (\nabla u_2-\nabla u_1\big )=0$ a.e. on $\Omega_1$. This means that $\nabla u_2=\nabla u_1$ a.e. on $\Omega_1$. So $\nabla (u_2-u_1)^{-}=\mathbf{0}$ a.e. on $\Omega$ which is a connected open set. Therefore $(u_2-u_1)^{-}$ is a constant positive function. If $(u_2-u_1)^{-}\equiv 0$ then $u_2-u_1=(u_2-u_1)^+\geq 0$ a.e. on $\Omega$ and we are done. If $(u_2-u_1)^{-}\equiv k$ for some constant $k>0$ then $u_2-u_1=-(u_2-u_1)^{-}=-k<0$ a.e. on $\Omega$ (so it is a strict negative function). But in this particular case we also have that $0=\displaystyle\int_{\Omega} c[(u_2-u_1)^-]^2\ dx=k^2\int_{\Omega} c\ dx\ \Rightarrow\ \int_{\Omega} c\ dx=0$ which implies that $c\equiv 0$ (since we knew that $c\geq 0$ a.e. on $\Omega$). Also $0=\displaystyle\int_{\Omega} (g_2-g_1)(u_2-u_1)^{-}\ dx=k\int_{\Omega} g_2-g_1\ dx\ \Rightarrow\ \int_{\Omega} g_2-g_1\ dx=0$. Since $g_2\geq g_1$ a.e. on $\Omega$ we finally obtain that $g_2\equiv g_1$ as needed.
	
	\noindent If $u_1\equiv u_2$ then relation \eqref{3ineqcomp} becomes:
	
	\begin{equation}
		0\geq\int_{\Omega} (g_2-g_1)\phi,\ \forall\ \phi\in W^{1,p(x)}(\Omega)^{+}.
	\end{equation}
	
	\noindent Since $\Omega$ is bounded we have that $\phi\equiv 1\in W^{1,p(x)}(\Omega)^{+}$. Hence we can write:
	
	\begin{equation}
		0\geq\int_{\Omega} \underbrace{g_2(x)-g_1(x)}_{\geq 0}\ dx\geq 0\ \Longrightarrow\ g_2\equiv g_1.
	\end{equation}
\end{proof}

\section{An auxiliary quasilinear elliptic problem}

\noindent In order to apply Rothe's method for the parabolic problem \eqref{eqpg} that will be studied in a future paper we need first to solve the following auxiliary quasilinear elliptic problem for any $\lambda\in (0,\infty)$ and $g\in \mathcal{M}_{\lambda}(\Omega)$:

\begin{equation}\tag{$E_{\lambda}$}\label{eqelambda}
	\begin{cases}-\operatorname{div}\mathbf{a}(x,\nabla V(x)) +\lambda V(x)=g(x), & x\in\Omega\\[3mm] \dfrac{\partial V}{\partial\nu}(x)=0, & x\in\partial\Omega\end{cases}
\end{equation}

\begin{definition} We say that $V\in W^{1,p(x)}(\Omega)$ is a \textbf{weak solution} for \eqref{eqelambda} if for any test function $\phi\in W^{1,p(x)}(\Omega)$ we have that:
	\begin{equation}
		\int_{\Omega} \mathbf{a}(x,\nabla V(x))\cdot \nabla\phi(x)\ dx+\lambda\int_{\Omega} V(x)\phi(x)\ dx=\int_{\Omega} g(x)\phi(x)\ dx.
	\end{equation}
	
\end{definition}

\begin{proposition} If $V\in W^{1,p(x)}(\Omega)$ is a weak solution of the problem \eqref{eqelambda} then $V\in\mathcal{U}$.
\end{proposition}

\begin{proof} Since $g\geq 0$ a.e. on $\Omega$ and $\lambda>0$ we deduce from the \textit{weak minimum principle} (Theorem \ref{weakmp}) that $V\geq 0$ a.e. on $\Omega$. Notice that for $1\in W^{1,p(x)}(\Omega)$ we have that:
	
	\begin{align*}
		\int_{\Omega} \mathbf{a}(x,\nabla 1)\cdot\nabla\phi\ dx+\lambda\int_{\Omega} 1\phi\ dx&=\int_{\Omega} \mathbf{a}(x,\nabla 1)\cdot\nabla\phi\ dx+\lambda\int_{\Omega} \phi\ dx=\int_{\Omega} \lambda\phi\ dx\\
		\big (g\in\mathcal{M}_{\lambda}(\Omega)\big ) \ \ \ \ \ \ \ 	&\geq\int_{\Omega} g(x)\phi(x)\ dx\\
		&=\int_{\Omega} \mathbf{a}(x,\nabla V(x))\cdot\nabla\phi(x)\ dx+\lambda\int_{\Omega} V(x)\phi(x)\ dx,\ \forall\ \phi\in W^{1,p(x)}(\Omega)^+.
	\end{align*}
	
	\noindent Thus from the \textit{weak comparison principle} we conclude that $1\geq V$ a.e. on $\Omega$. In conclusion $V\in W^{1,p(x)}\cap\mathcal{U}$.
\end{proof}

\section{The perturbed problem}
\noindent In order to show the existence for \eqref{eqelambda} we will use some information about the perturbed problem:

\begin{equation}\tag{$E_{\lambda,\varepsilon}$}\label{eqelambdaeps}
	\begin{cases}-\operatorname{div}\mathbf{a}(x,\nabla V(x)) +\lambda V(x)+\varepsilon |V(x)|^{p(x)-2}V(x)=g(x), & x\in\Omega\\[3mm] \dfrac{\partial V}{\partial\nu}(x)=0, & x\in\partial\Omega\end{cases}
\end{equation}

\noindent where $\varepsilon\geq 0$. Note that for $\varepsilon=0$ problem \eqref{eqelambdaeps} is exactly problem \eqref{eqelambda}.

\begin{definition} We say that $V_{\varepsilon}\in W^{1,p(x)}(\Omega)$ is a \textbf{weak solution} for \eqref{eqelambdaeps} if for any test function $\phi\in W^{1,p(x)}(\Omega)$ we have that:
	\begin{equation}
		\int_{\Omega} \mathbf{a}(x,\nabla V_{\varepsilon}(x))\cdot \nabla\phi(x)\ dx+\lambda\int_{\Omega} V_{\varepsilon}(x)\phi(x)\ dx+\varepsilon\int_{\Omega} |V_{\varepsilon}(x)|^{p(x)-2}V_{\varepsilon}(x)\phi(x)\ dx=\int_{\Omega} g(x)\phi(x)\ dx.
	\end{equation}
	
\end{definition}

\begin{remark} We know that $V_{\varepsilon}\in W^{1,p(x)}(\Omega)\subset L^{p(x)}(\Omega)$. Therefore $|V_{\varepsilon}|^{p(x)-1}\in L^{\frac{p(x)}{p(x)-1}}(\Omega)=L^{p'(x)}(\Omega)$.\footnote{See Theorem \ref{athleb} \textbf{(24)} from the Appendix.} Now since $\phi\in W^{1,p(x)}(\Omega)\subset L^{p(x)}(\Omega)$, from \textit{H\"{o}lder inequality} we get that:
	
	\begin{equation*}
		\left |\int_{\Omega} |V_{\varepsilon}(x)|^{p(x)-2}V_{\varepsilon}(x)\phi(x)\ dx \right |\leq \int_{\Omega} |V_{\varepsilon}|^{p(x)-1}\phi(x)\ dx\leq 2\Vert V_{\varepsilon}^{p(x)-1}\Vert_{L^{p'(x)}(\Omega)}\Vert\phi\Vert_{L^{p(x)}(\Omega)}<\infty.
	\end{equation*}
	\noindent So, $\Omega\ni x\longmapsto |V_{\varepsilon}(x)|^{p(x)-2}V_{\varepsilon}(x)\phi(x)$ is in $L^1(\Omega)$. 
	
\end{remark}

\begin{proposition} If $V_{\varepsilon}\in W^{1,p(x)}(\Omega)$ is a weak solution of the problem \eqref{eqelambdaeps} then $V_{\varepsilon}\in\mathcal{U}$.
\end{proposition}

\begin{proof} First we show that $V_{\varepsilon}\geq 0$ a.e. on $\Omega$. Take $\phi=V_{\varepsilon}^{-}\in W^{1,p(x)}(\Omega)^{+}$ as the test function. Then $\nabla \phi=\nabla V_{\varepsilon}^{-}=\begin{cases}0,\ \text{a.e. on}\ \{x\in\Omega\ |\ V_{\varepsilon}(x)\geq 0\}\\[3mm]-\nabla V_{\varepsilon}, \ \text{a.e. on}\ \{x\in\Omega\ |\ V_{\varepsilon}(x)<0\} \end{cases}$. We obtain that:
	
	\begin{equation}0\geq -\int_{\Omega} \Phi(x,|\nabla V_{\varepsilon}|)|\nabla V_{\varepsilon}^-|\ dx-\lambda\int_{\Omega} (V_{\varepsilon}^{-})^2\ dx-\varepsilon\int_{\Omega} |V_{\varepsilon}|^{p(x)-2}(V_{\varepsilon}^{-})^2\ dx=\int_{\Omega} gV_{\varepsilon}^{-}\ dx\geq 0
	\end{equation}
	
	\noindent As $\lambda>0$ we get that $\displaystyle\int_{\Omega} (V_{\varepsilon}^{-})^2\ dx$, i.e. $V_{\varepsilon}^{-}\equiv 0$. So $V_{\varepsilon}=V_{\varepsilon}^+\geq 0$ a.e. on $\Omega$.
	
	\noindent For the second part of the proof first note that $\forall\ \phi\in W^{1,p(x)}(\Omega)$:
	
	\begin{equation}
		\begin{cases}\displaystyle\int_{\Omega}\mathbf{a}(x,\nabla 1)\cdot\nabla\phi(x)\ dx+\lambda\int_{\Omega} 1\cdot\phi(x)\ dx+\varepsilon\int_{\Omega} 1^{p(x)-1}\phi(x)\ dx=\int_{\Omega} (\lambda+\varepsilon)\phi(x)\ dx\\[5mm]
			\displaystyle\int_{\Omega} \mathbf{a}(x,\nabla V_{\varepsilon}(x))\cdot\nabla\phi(x)\ dx+\lambda\int_{\Omega} V_{\varepsilon}(x)\phi(x)\ dx+\varepsilon\int_{\Omega} |V_{\varepsilon}(x)|^{p(x)-2}V_{\varepsilon}(x)\phi(x)\ dx=\int_{\Omega} g(x)V_{\varepsilon}(x)\ dx\end{cases}.
	\end{equation}
	
	\noindent Substracting both equations, and using the fact that $V_{\varepsilon}\geq 0$ a.e. on $\Omega$ we have that: $|V_{\varepsilon}|^{p(x)-2}V_{\varepsilon}=V_{\varepsilon}^{p(x)-1}$, and this leads to:
	
	\begin{equation}
		\int_{\Omega} \big (\mathbf{a}(x,\nabla 1)-\mathbf{a}(x,\nabla V_{\varepsilon})\big )\cdot\nabla\phi\ dx+\lambda\int_{\Omega}\big (1-V_{\varepsilon}\big )\phi\ dx+\varepsilon\int_{\Omega} \big (1-V_{\varepsilon}^{p(x)-1}\big )\phi\ dx=\int_{\Omega} (\lambda+\varepsilon-g)\phi\ dx.
	\end{equation}

	\noindent We choose $\phi=(1-V_{\varepsilon})^{-}=\begin{cases} 0, \ \text{a.e. on}\ \{x\in\Omega\ |\ V_{\varepsilon}(x)\leq 1\}\\[3mm] -(1-V_{\varepsilon}(x)),\ \text{a.e. on}\ \{x\in\Omega\ |\ V_{\varepsilon}(x)>1\}:=\Omega_1\end{cases}\in W^{1,p(x)}(\Omega)$ as the test function. Notice that:
	
	\begin{equation}
		\nabla\phi=\nabla (1-V_{\varepsilon})^{-}=\begin{cases} 0,\ \text{a.e. on}\ \{x\in\Omega\ |\ V_{\varepsilon}(x)\leq 1\}\\[3mm] -(\nabla 1-\nabla V_{\varepsilon})=\nabla V_{\varepsilon},\ \text{a.e. on}\ \Omega_1\end{cases}.
	\end{equation}
	
	\noindent We obtain that:
	
	\begin{equation}
		-\int_{\Omega_1} \underbrace{\mathbf{a}(x,\nabla V_{\varepsilon})\cdot\nabla V_{\varepsilon}\ dx}_{\geq 0}-\lambda\int_{\Omega} [(1-V_{\varepsilon})^{-}]^2\ dx-\varepsilon\int_{\Omega_1} \big (1-V_{\varepsilon}^{p(x)-1} \big )\cdot (1-V_{\varepsilon})\ dx=\int_{\Omega} (\lambda+\varepsilon-g) (1-V_{\varepsilon})^{-}\ dx.
	\end{equation}

	\noindent It is straightforward to check that for each $x\in\overline{\Omega}$ the function $\theta\mapsto \theta^{p(x)-1}$ is a strictly increasing function on $[0,\infty)$, as $p(x)\geq p^->1$ for all $x\in\overline{\Omega}$. So the left hand side is negative while the right hand side is positive as $\lambda+\varepsilon\geq\lambda\geq g$. Therefore $\lambda\displaystyle\int_{\Omega} [(1-V_{\varepsilon})^{-}]^2\ dx=0$ which shows that $(1-V_{\varepsilon})^{-}\equiv 0$ on $\Omega$ ($\lambda>0$). In conclusion $1-V_{\varepsilon}=(1-V_{\varepsilon})^{+}\geq 0$ a.e. on $\Omega$.
\end{proof}

\begin{theorem}\label{3thsoleps}
	Problem \eqref{eqelambdaeps} admits a unique weak solution for every $g\in\mathcal{M}_{\lambda}(\Omega)$ that will be denoted $V_{\varepsilon}\in\mathcal{U}$ for each $\varepsilon>0$.
\end{theorem}

\begin{proof} \textbf{(Uniqueness)} Suppose that $V_1,V_2\in W^{1,p(x)}(\Omega)$ are two weak solutions of \eqref{eqelambda}. Choosing $\phi=V_1-V_2\in W^{1,p(x)}(\Omega)$ as the test function we easily get that:
	
	\begin{equation}
		\int_{\Omega} \big (\mathbf{a}(x,\nabla V_1)-\mathbf{a}(x,\nabla V_2) \big)\cdot (\nabla V_1-\nabla V_2)\ dx+\lambda\int_{\Omega} (V_1-V_2)^2\ dx+\varepsilon\int_{\Omega} \big (|V_1|^{p(x)-2}V_1-|V_2|^{p(x)-2}V_2\big )(V_1-V_2)\ dx=0.
	\end{equation}
	
	\noindent Using the fact that $\mathbf{a}(x,\cdot)$ is monotone (see Proposition \ref{4prop1}), $\lambda>0$ and the function $\mathbb{R}\ni\theta\mapsto |\theta|^{p(x)-2}\theta=|\theta|^{p(x)-1}\operatorname{sgn}(\theta)$ is strictly increasing on $\mathbb{R}$ for all $x\in\overline{\Omega}$ as $p(x)\geq p^{-}>1,\ \forall x\in\overline{\Omega}$, we obtain that every integral in the left hand side is positive. We infer that $\displaystyle\int_{\Omega} (V_1-V_2)^2\ dx=0$. So $V_1\equiv V_2$, as needed.

	\noindent\textbf{(Existence)} Consider $\mathcal{J}_{\lambda,\varepsilon}:W^{1,p(x)}(\Omega)\to\mathbb{R}$ the Euler-Lagrange (or energy) functional associated to \eqref{eqelambda}, which is given by:
	
	\begin{equation}
		\mathcal{J}_{\lambda,\varepsilon}(V)=\mathcal{A}(V)+\lambda\mathcal{I}(V)+\varepsilon\mathcal{C}(V)-\mathcal{G}(V),
	\end{equation}
	\noindent where:
	
	\begin{equation}
		\mathcal{I},\mathcal{C},\mathcal{G}:W^{1,p(x)}(\Omega)\to\mathbb{R},\ \begin{cases} \mathcal{I}(V)=\dfrac{1}{2}\displaystyle\int_{\Omega} V^2(x)\ dx,\ \mathcal{C}(V)=\int_{\Omega}\dfrac{|V(x)|^{p(x)}}{p(x)}\ dx \\[5mm] \mathcal{G}(V)=\displaystyle\int_{\Omega} g(x)V(x)\ dx\end{cases}.
	\end{equation}	
	
	\noindent So, to see it better:
	
	\begin{equation}
		\mathcal{J}_{\lambda,\varepsilon}(V)=\int_{\Omega} A(x,\nabla V(x))\ dx+\dfrac{\lambda}{2}\int_{\Omega} V^2(x)\ dx+\varepsilon\int_{\Omega}\dfrac{|V(x)|^{p(x)}}{p(x)}\ dx-\int_{\Omega}g(x)V(x)\ dx,\ \forall\ V\in W^{1,p(x)}(\Omega).
	\end{equation}
	
	\noindent $\blacktriangleright$ \textbf{Fact I: $\mathcal{J}_{\lambda,\varepsilon}$ is well-defined and bounded from below.}
	
	\noindent We have prove before (see Proposition \ref{4prop1} \textbf{(6)}) that $\mathcal{A}$ is well-defined. Because $2<p^*(x)$ for all $x\in\overline{\Omega}$ we get that $W^{1,p(x)}(\Omega)\hookrightarrow L^2(\Omega)$\footnote{See Theorem \ref{athmemb} from the Appendix.}, so $0\leq \mathcal{I}(V)=\dfrac{\lambda}{2}\displaystyle\int_{\Omega} V^2(x)\ dx=\dfrac{\lambda}{2}\Vert V\Vert_{L^2(\Omega)}\lesssim \Vert V\Vert_{W^{1,p(x)}(\Omega)} <\infty,\ \forall\ V\in W^{1,p(x)}(\Omega)$, which shows that $\mathcal{I}$ is well-defined. 
	
	\noindent Going further, from $W^{1,p(x)}(\Omega)\subset L^{p(x)}(\Omega)$, we get for each $V\in W^{1,p(x)}(\Omega)$ that: $0\leq \mathcal{C}(V)=\displaystyle\int_{\Omega} \dfrac{|V(x)|^{p(x)}}{p(x)}\ dx\leq\dfrac{1}{p^{-}}\int_{\Omega} |V(x)|^{p(x)}\ dx=\dfrac{1}{p^{-}}\rho_{p(x)}(V)<\infty$, from the definition of the space $L^{p(x)}(\Omega)$. So $\mathcal{C}$ is well-defined.
	
	\noindent To see that the functional $\mathcal{G}$ is well-defined we say that:
	
	\begin{align*}
		|\mathcal{G}(V)|&=\left |\int_{\Omega} g(x)V(x)\ dx \right |\leq \int_{\Omega} |g(x)|\cdot |V(x)|\ dx\\
		\big (g\in\mathcal{M}_{\lambda}(\Omega)\big )\ \ \ \ \ \ \ 	&\leq \lambda\Vert V\Vert_{L^1(\Omega)}\lesssim \Vert V\Vert_{W^{1,p(x)}(\Omega)}<\infty,\ \forall\ V\in W^{1,p(x)}(\Omega),
	\end{align*}
	
	\noindent as $W^{1,p(x)}(\Omega)\hookrightarrow L^1(\Omega)$. Remark also that $\mathcal{G}\in W^{1,p(x)}(\Omega)^*$ i.e. it is a bounded linear operator. Hence $\mathcal{J}_{\lambda,\varepsilon}$ is well-defined.
	
	\noindent At last $\mathcal{J}_{\lambda,\varepsilon}$ is bounded from below because:
	
	\begin{align*}
		\mathcal{J}_{\lambda,\varepsilon}(V)&=\int_{\Omega}A(x,\nabla V)\ dx+\dfrac{\lambda}{2}\int_{\Omega} V^2\ dx+\varepsilon\int_{\Omega}\dfrac{|V|^{p(x)}}{p(x)}\ dx-\int_{\Omega}gV\ dx\\
		\text{(Cauchy ineq.)}	\ \ \ \ \ \ \ \	&\geq 0+\dfrac{\lambda}{2}\Vert V\Vert^2_{L^2(\Omega)}+0-\Vert g\Vert_{L^2(\Omega)}\Vert V\Vert_{L^2(\Omega)}\\
		&=\dfrac{\lambda}{2}\left (\Vert V\Vert_{L^2(\Omega)}-\dfrac{\Vert g\Vert_{L^2(\Omega)}}{\lambda} \right)^2-\dfrac{\Vert g\Vert^2_{L^2(\Omega)}}{2\lambda}\\
		&\geq -\dfrac{\Vert g\Vert^2_{L^2(\Omega)}}{2\lambda}>-\infty,\ \forall\ V\in W^{1,p(x)}(\Omega).
	\end{align*}
	
	\bigskip
	
	\noindent$\blacktriangleright$ \textbf{Fact II:} $\mathcal{J}_{\lambda,\varepsilon}\in C^1\big (W^{1,p(x)}(\Omega)\big )$ and $\forall\ \phi\in W^{1,p(x)}(\Omega)$ the following formula holds:
	
	\begin{equation}\label{3eqc}
		\langle\mathcal{J}'_{\lambda,\varepsilon}(V),\phi\rangle=\int_{\Omega}\mathbf{a}(x,\nabla V)\cdot\phi(x)\ dx+\lambda\int_{\Omega} V\phi\ dx+\varepsilon\int_{\Omega} |V(x)|^{p(x)-2}V\phi\ dx-\int_{\Omega} gV\ dx.
	\end{equation}
	
	\begin{remark} From \eqref{3eqc} we deduce that $V_{\varepsilon}$ is the weak solution of the problem \eqref{eqelambdaeps} if and only if $\mathcal{J}'_{\lambda,\varepsilon}(V_{\varepsilon})=0$, i.e. $V_{\varepsilon}$ is a critical point of the functional $\mathcal{J}_{\lambda,\varepsilon}\in C^1\big (W^{1,p(x)}(\Omega) \big )$.
	\end{remark}
	
	\noindent We already know from Proposition \ref{4prop1} \textbf{(6)} that $\mathcal{A}\in C^1\big (W^{1,p(x)}(\Omega)\big )$ and:
	
	\begin{equation}\label{3eqc1}
		\langle\mathcal{A}'(V),\phi\rangle=\int_{\Omega}\mathbf{a}(x,\nabla V(x))\cdot\nabla\phi(x)\ dx,\ \forall\ V, \phi\in W^{1,p(x)}(\Omega).
	\end{equation}
	
	\noindent Fix some $V\in W^{1,p(x)}(\Omega)$ and consider the linear operator $L_V:W^{1,p(x)}(\Omega)\to\mathbb{R},\ L_V(\phi)=\displaystyle\int_{\Omega} V(x)\phi(x)\ dx$.  First we show that it is bounded. From the fact that $W^{1,p(x)}(\Omega)\hookrightarrow L^2(\Omega)$ we deduce that $\exists\ \tilde{C}>0$ such that $\Vert U\Vert_{L^2(\Omega)}\leq \tilde{C}\Vert U\Vert_{W^{1,p(x)}(\Omega)},\ \forall\ U\in W^{1,p(x)}(\Omega)$. Hence:

	\begin{align*}
		|L_{V}(\phi)|&=\left |\int_{\Omega}V(x)\phi(x)\ dx\right |\leq \int_{\Omega} |V(x)|\cdot |\phi(x)|\ dx\\
		\text{(Cauchy ineq.)}\ \ \ \ \ \ \ 	&\leq \Vert V\Vert_{L^2(\Omega)}\cdot\Vert\phi\Vert_{L^2(\Omega)}\\
		\big (W^{1,p(x)}(\Omega)\hookrightarrow L^2(\Omega) \big )\ \ \ \ \ \ \ \ &\leq \tilde{C}^2\Vert V\Vert_{W^{1,p(x)}(\Omega)}\cdot \Vert \phi\Vert_{W^{1,p(x)}(\Omega)}
	\end{align*}
	
	\noindent So $\Vert L_V\Vert_{W^{1,p(x)}(\Omega)^*}=\displaystyle\sup_{\phi\in W^{1,p(x)}(\Omega)\setminus\{0\}}\dfrac{|L_V(\phi)|}{\Vert\phi\Vert_{W^{1,p(x)}(\Omega)}}\leq\tilde{C}^2\Vert V\Vert_{W^{1,p(x)}(\Omega)}<\infty$, which shows that $L_V\in W^{1,p(x)}(\Omega)^*$ for each $V\in W^{1,p(x)}(\Omega)$. Moreover:
	
	\begin{align*}
		\dfrac{\big |\mathcal{I}(V+\phi)-\mathcal{I}(V)-L_V(\phi)\big |}{\Vert\phi\Vert_{W^{1,p(x)}(\Omega)}}=\dfrac{\left |\displaystyle\int_{\Omega} \phi^2(x)\ dx \right |}{\Vert\phi\Vert_{W^{1,p(x)}(\Omega)}}=\dfrac{\Vert\phi\Vert^2_{L^{2}(\Omega)}}{\Vert\phi\Vert_{W^{1,p(x)}(\Omega)}}\leq \tilde{C}^2\dfrac{\Vert\phi\Vert^2_{W^{1,p(x)}(\Omega)}}{\Vert\phi\Vert_{W^{1,p(x)}(\Omega)}}=\tilde{C}^2\Vert\phi\Vert_{W^{1,p(x)}(\Omega)}\longrightarrow\ 0,
	\end{align*}
	
	\noindent as $\phi\to 0$ in $W^{1,p(x)}(\Omega)$. This shows that $\mathcal{I}$ is Fr\'{e}chet differentiable and $\mathcal{I}'(V)=L_V$. 
	
	\noindent The last thing to show here is that $\mathcal{I}':W^{1,p(x)}(\Omega)\to W^{1,p(x)}(\Omega)^*$ is a continuous operator. In that sense consider $(V_n)_{n\geq 1}\subset W^{1,p(x)}(\Omega)$ with $\Vert V_n-V\Vert_{W^{1,p(x)}(\Omega)}\to 0$. We have that:
	
	\begin{align*}
		\Vert\mathcal{I}(V_n)-\mathcal{I}(V)\Vert_{W^{1,p(x)}(\Omega)^*}&=\sup_{\phi\in W^{1,p(x)}(\Omega)\setminus\{0\}}\dfrac{\left |\displaystyle\int_{\Omega} (V_n(x)-V(x))\phi(x)\ dx \right |}{\Vert\phi\Vert_{W^{1,p(x)}(\Omega)}}\\
		\text{(Cauchy ineq.)}\ \ \ \ \ \ \ 		&\leq \sup_{\phi\in W^{1,p(x)}(\Omega)\setminus\{0\}}\dfrac{\Vert V_n-V\Vert_{L^2(\Omega)}\Vert \phi\Vert_{L^2(\Omega)}}{\Vert\phi\Vert_{W^{1,p(x)}(\Omega)}}\\
		\big (W^{1,p(x)}(\Omega)\hookrightarrow L^2(\Omega)  \big )\ \ \ \ \ \ \ &\leq \tilde{C}^2\sup_{\phi\in W^{1,p(x)}(\Omega)\setminus\{0\}}\dfrac{\Vert V_n-V\Vert_{W^{1,p(x)}(\Omega)}\Vert \phi\Vert_{W^{1,p(x)}(\Omega)}}{\Vert\phi\Vert_{W^{1,p(x)}(\Omega)}}\\
		&=\tilde{C}^2\Vert V_n-V\Vert_{W^{1,p(x)}(\Omega)}\stackrel{n\to\infty}{\longrightarrow} 0.
	\end{align*}
	
	\noindent In conclusion $\mathcal{I}'\in C^1\big (W^{1,p(x)}(\Omega) \big )$ and:
	
	\begin{equation}\label{3eqc2}
		\langle \mathcal{I}'(V),\phi\rangle=\int_{\Omega} V(x)\phi(x)\ dx,\ \forall\ V,\phi\in W^{1,p(x)}(\Omega).
	\end{equation}

	\noindent For a fixed $V\in W^{1,p(x)}(\Omega)$ we consider the linear operator $\ell_V:W^{1,p(x)}(\Omega)\to\mathbb{R}$ given by:
	
	\begin{equation}
		\ell_V(\phi)=\int_{\Omega} |V(x)|^{p(x)-2}V(x)\phi(x)\ dx,\ \forall\ \phi\in W^{1,p(x)}(\Omega).
	\end{equation}
	
	\noindent The first think to observe is that $\ell_V\in W^{1,p(x)}(\Omega)^*$. Indeed, for all $\phi\in W^{1,p(x)}(\Omega)$:
	
	\begin{equation*}
		\left | \ell_V(\phi) \right |\leq \int_{\Omega} |V|^{p(x)-1}\phi(x)\ dx\leq 2\Vert V^{p(x)-1}\Vert_{L^{p'(x)}(\Omega)}\Vert\phi\Vert_{L^{p(x)}(\Omega)}\leq 2\Vert V^{p(x)-1}\Vert_{L^{p'(x)}(\Omega)}\Vert\phi\Vert_{W^{1,p(x)}(\Omega)}.
	\end{equation*}
	
	\noindent So $\displaystyle\sup_{\phi\in W^{1,p(x)}(\Omega)\setminus\{0\}}\dfrac{|\ell_V(\phi)|}{\Vert\phi\Vert_{W^{1,p(x)}(\Omega)}}\leq 2\Vert V^{p(x)-1}\Vert_{L^{p'(x)}(\Omega)}\ \Rightarrow\ \ell_V\in W^{1,p(x)}(\Omega)^*$.

	\noindent Fix some $V\in W^{1,p(x)}(\Omega)$ and a direction $\phi\in W^{1,p(x)}(\Omega)$. For any sequence of real non-zero numbers $(\varepsilon_n)_{n\geq 1}$ convergent to $0$ we define the sequence of function $(v_n)_{n\geq 1}$ by:
	
	\begin{equation}
		v_n(x)=\dfrac{|V(x)+\varepsilon_n\phi(x)|^{p(x)}-|V(x)|^{p(x)}}{p(x)\varepsilon_n},\ \text{for a.a.}\ x\in\Omega,\ \forall\ n\geq 1.
	\end{equation}
	
	\noindent Using Lemma \ref{alemhos} ($p(x)\geq p^->1$ for all $x\in\overline{\Omega}$) we have that: $\lim\limits_{n\to\infty} v_n(x)=|V(x)|^{p(x)-2}V(x)\phi(x)$ for a.a. $x\in\Omega$.

	\noindent Taking into account that both $V+\varepsilon_n\phi, \phi\in W^{1,p(x)}(\Omega)\subset L^{p(x)}(\Omega)$ we deduce from Theorem \ref{athleb} \textbf{(24)} that $|V+\varepsilon\phi|^{p(x)},|V(x)|^{p(x)}\in L^1(\Omega)$. Thence $v_n\in L^1(\Omega)$ for each integer $n\geq 1$.
	
	\noindent From the \textit{Lebesgue dominated convergence theorem} we obtain that:
	
	\begin{align*}
		\lim\limits_{n\to\infty}\dfrac{\mathcal{C}(V+\varepsilon_n\phi)-\mathcal{C}(V)}{\varepsilon_n}&=\lim\limits_{n\to\infty}\int_{\Omega}\dfrac{|V(x)+\varepsilon_n\phi(x)|^{p(x)}-|V(x)|^{p(x)}}{p(x)\varepsilon_n}\ dx\\
		&=\lim\limits_{n\to\infty}\int_{\Omega} v_n(x)\ dx=\int_{\Omega} |V(x)|^{p(x)-2}V(x)\phi(x)\ dx\\
		&=\ell_V(\phi),\ \forall\ V,\phi\in W^{1,p(x)}(\Omega).
	\end{align*}
	
	\noindent So $\partial\mathcal{C}(V):W^{1,p(x)}(\Omega)\to\mathbb{R},\ \partial\mathcal{C}(V)\phi=\ell_{V}(\phi)$ is a bounded linear operator as we have seen before. Thus we can define the G\^ateaux differential of $\mathcal{C}$ as $\partial\mathcal{C}:W^{1,p(x)}(\Omega)\to W^{1,p(x)}(\Omega)^*$.
	
	\noindent Next, we need to show that $\partial\mathcal{C}$ is a continuous (nonlinear) operator. In that sense consider some $V\in W^{1,p(x)}(\Omega)$ and any sequence $(V_n)_{n\geq 1}\subset W^{1,p(x)}(\Omega)$ with $V_n\longrightarrow V$ in $W^{1,p(x)}(\Omega)$. Also denote: $\Omega_1=\{x\in\overline{\Omega}\ |\ p(x)<2\}$ and $\Omega_2=\{x\in\overline{\Omega}\ |\ p(x)\geq 2\}$. Therefore:
	
	\begin{align*}
		&\Vert \partial \mathcal{C}(V_n)-\partial\mathcal{C}(V)\Vert_{W^{1,p(x)}(\Omega)^*}=\sup_{\phi\in W^{1,p(x)}(\Omega)\setminus\{0\}}\dfrac{\left |\displaystyle\int_{\Omega} \big (|V_n(x)|^{p(x)-2}V_n(x)-|V(x)|^{p(x)-2}V(x)\big )\phi(x) \ dx\right |}{\Vert\phi\Vert_{W^{1,p(x)}(\Omega)}}\\
		&\leq \sup_{\phi\in W^{1,p(x)}(\Omega)\setminus\{0\}}\dfrac{\displaystyle\int_{\Omega} \left | |V_n(x)|^{p(x)-2}V_n(x)-|V(x)|^{p(x)-2}V(x)\right |\cdot |\phi(x)| \ dx}{\Vert\phi\Vert_{W^{1,p(x)}(\Omega)}}\\
		&\leq\sup_{\phi\in W^{1,p(x)}(\Omega)\setminus\{0\}}\dfrac{\displaystyle\int_{\Omega_1} \left | |V_n(x)|^{p(x)-2}V_n(x)-|V(x)|^{p(x)-2}V(x)\right |\cdot |\phi(x)| \ dx}{\Vert\phi\Vert_{W^{1,p(x)}(\Omega)}}+\\
		&+\dfrac{\displaystyle\int_{\Omega_2} \left | |V_n(x)|^{p(x)-2}V_n(x)-|V(x)|^{p(x)-2}V(x)\right |\cdot |\phi(x)| \ dx}{\Vert\phi\Vert_{W^{1,p(x)}(\Omega)}}\\
		\text{(Lemma \ref{alplap})}		&\leq\sup_{\phi\in W^{1,p(x)}(\Omega)\setminus\{0\}}\dfrac{c_1\displaystyle\int_{\Omega_1} |V_n-V|^{p(x)-1}|\phi(x)|\ dx}{\Vert\phi\Vert_{W^{1,p(x)}(\Omega)}}+\dfrac{c_2\displaystyle\int_{\Omega_2}\big (|V_n|+|V| \big )^{p(x)-2} |V_n-V|\cdot|\phi(x)|\ dx}{\Vert\phi\Vert_{W^{1,p(x)}(\Omega)}}\\
	\end{align*}
	
	\noindent We have here two separated problems. From $|V_n-V|\in W^{1,p(x)}(\Omega)\subset L^{p(x)}(\Omega)$ we deduce from Theorem \ref{athleb} \textbf{(24)} that $|V_n-V|^{p(x)-1}\in L^{p'(x)}(\Omega)$ for each $n\geq 1$. Since $|\phi|\in W^{1,p(x)}(\Omega)\subset L^{p(x)}(\Omega)$ we get from \textit{H\"{o}lder inequality} that:
	
	\begin{align*}
		\dfrac{1}{\Vert\phi\Vert_{W^{1,p(x)}(\Omega)}}\int_{\Omega_1} |V_n-V|^{p(x)-1}|\phi(x)|\ dx &\leq \dfrac{2}{\Vert\phi\Vert_{W^{1,p(x)}(\Omega)}}\Vert |V_n-V|^{p(x)-1}\Vert_{L^{p'(x)}(\Omega_1)}\Vert\phi\Vert_{L^{p(x)}(\Omega_1)}\\
		&\leq \dfrac{2}{\Vert\phi\Vert_{W^{1,p(x)}(\Omega)}}\Vert |V_n-V|^{p(x)-1}\Vert_{L^{p'(x)}(\Omega)}\Vert\phi\Vert_{L^{p(x)}(\Omega)}\\
		&\leq  2\Vert |V_n-V|^{p(x)-1}\Vert_{L^{p'(x)}(\Omega)}\stackrel{n\to\infty}{\longrightarrow} 0\\
		\text{(Theorem 2.58 from \cite{Cruz})}\ \ \ \ \ \ \ 	&\Longleftrightarrow\rho_{p'(x)}\left (|V_n-V|^{p(x)-1} \right )\stackrel{n\to\infty}{\longrightarrow} 0\\
		&\Longleftrightarrow\int_{\Omega}\left (|V_n-V|^{p(x)-1} \right )^{p'(x)}\ dx\stackrel{n\to\infty}{\longrightarrow} 0\\
		&\Longleftrightarrow\int_{\Omega}|V_n-V|^{p(x)}\ dx\stackrel{n\to\infty}{\longrightarrow} 0\\
		&\Longleftrightarrow\rho_{p(x)}\left (V_n-V \right )\stackrel{n\to\infty}{\longrightarrow} 0\\
		\text{(Theorem 2.58 from \cite{Cruz})}\ \ \ \ \ \ \	&\Longleftrightarrow\Vert V_n-V\Vert_{L^{p(x)}(\Omega)}\stackrel{n\to\infty}{\longrightarrow} 0,
	\end{align*}
	
	\noindent which is true, because $V_n\to V$ in $W^{1,p(x)}(\Omega)$ so that: $0\leq \Vert V_n-V\Vert_{L^{p(x)}(\Omega)}\leq\Vert V_n-V\Vert_{W^{1,p(x)}(\Omega)} \stackrel{n\to\infty}{\longrightarrow} 0$. 
	
	\noindent For the other term notice that $z_n:=|V_n|+|V|\in W^{1,p(x)}(\Omega)\subset L^{p(x)}(\Omega)$, and therefore, from Theorem \ref{athleb} \textbf{(24)}, we have that $z_n^{p(x)-2}\in L^{\frac{p(x)}{p(x)-2}}(\Omega)$ for each $n\geq 1$. Another important remark is that $z_n\to 2|V|$ in $L^{p(x)}(\Omega)$ because $\big |z_n-2|V|\big |=\big ||V_n|-|V|\big |\leq |V_n-V|$, which implies that $0\leq\rho_{p(x)}\big (z_n-2|V| \big )\leq \rho_{p(x)}\big (V_n-V \big )\stackrel{n\to\infty}{\longrightarrow} 0$.\footnote{See Theorem \ref{athleb} \textbf{(14)}.}
	
	\noindent From the \textit{Generalized H\"{o}lder inequality}\footnote{It can be found at \cite[Corollary 2.30, page 30]{Cruz}.}, taking into account that $\left (\dfrac{p(x)}{p(x)-2} \right )^{-1}+\dfrac{1}{p(x)}+\dfrac{1}{p(x)}=1$ we obtain that:
	
	\begin{align*}
		\dfrac{1}{\Vert\phi\Vert_{W^{1,p(x)}(\Omega)}}\int_{\Omega_2}\big (|V_n|+|V| \big )^{p(x)-2} |V_n-V|\cdot|\phi(x)|\ dx &\lesssim \Vert z_n^{p(x)-2}\Vert_{L^{\frac{p(x)}{p(x)-2}}(\Omega_2)}\Vert V_n-V\Vert_{L^{p(x)}(\Omega_2)}\dfrac{\Vert\phi\Vert_{L^{p(x)}(\Omega_2)}}{\Vert\phi\Vert_{W^{1,p(x)}(\Omega)}}      \\
		&\leq\Vert z_n^{p(x)-2}\Vert_{L^{\frac{p(x)}{p(x)-2}}(\Omega)}\Vert V_n-V\Vert_{L^{p(x)}(\Omega)}\dfrac{\Vert\phi\Vert_{L^{p(x)}(\Omega)}}{\Vert\phi\Vert_{W^{1,p(x)}(\Omega)}} \\
		&\leq \Vert z_n^{p(x)-2}\Vert_{L^{\frac{p(x)}{p(x)-2}}(\Omega)}\Vert V_n-V\Vert_{L^{p(x)}(\Omega)}\stackrel{n\to\infty}{\longrightarrow} 0,
	\end{align*}
	
	\noindent because the sequence of real positive numbers $\left (\Vert z_n^{p(x)-2}\Vert_{L^{\frac{p(x)}{p(x)-2}}(\Omega)} \right )_{n\geq 1}$ is bounded. Let's explain this fact in details: From $z_n\to 2|V|$ in $L^{p(x)}(\Omega)$ we get that there is some $n_0\geq 1$ such that:
	
	\begin{equation}
		\Vert z_n\Vert_{L^{p(x)}(\Omega)}\leq 2\Vert V\Vert_{L^{p(x)}(\Omega)}+\Vert z_n-2|V|\Vert_{L^{p(x)}(\Omega)}\leq 2\Vert V\Vert_{L^{p(x)}(\Omega)}+1:=\alpha,\ \forall\ n\geq n_0.
	\end{equation}
	
	\noindent Thus $\left\Vert\dfrac{z_n}{\alpha}\right \Vert_{L^{p(x)}(\Omega)}\leq 1$ for any $n\geq n_0$. Using now Theorem \ref{athleb} \textbf{(25)} we deduce that $\rho_{p(x)}\left (\dfrac{z_n}{\alpha} \right )\leq \left\Vert\dfrac{z_n}{\alpha}\right \Vert_{L^{p(x)}(\Omega)}\leq 1$ for $n\geq n_0$. Because we have chosen $\alpha\geq 1$ we observe that:
	
	\begin{equation}
		\int_{\Omega} \dfrac{z_n^{p(x)}}{\alpha^{p^+}}\ dx\leq \int_{\Omega}\left (\dfrac{z_n}{\alpha}\right )^{p(x)}\ dx=\rho_{p(x)}\left (\dfrac{z_n}{\alpha} \right )\leq 1,\ \forall\ n\geq n_0.
	\end{equation}
	
	\noindent On the other hand 
	
	\begin{align*}
		\Vert z_n^{p(x)-2}\Vert_{L^{\frac{p(x)}{p(x)-2}}(\Omega)}&=\inf\left\{ \lambda>0\ |\ \rho_{\frac{p(x)}{p(x)-2}}\left(z_n^{p(x)-2}\right)\leq 1\right\}\\
		&=\inf\underbrace{\left\{ \lambda>0\ \bigg |\ \int_{\Omega}\dfrac{z_n^{p(x)}}{\lambda^{\frac{p(x)}{p(x)-2}}}\ dx\leq 1 \right\}}_{:=Z_n}
	\end{align*}
	
	\noindent If we choose $\lambda>1$ sufficiently large so that $\lambda^{\frac{p^+}{p^+-2}}\geq\alpha^{p^+}$, i.e. $\lambda\geq\alpha^{p^+-2}$ then:
	
	\begin{equation}
		\lambda^{\frac{p(x)}{p(x)-2}}=\lambda^{1+\frac{2}{p(x)-2}}\stackrel{\lambda>1}{\geq}\lambda^{1+\frac{2}{p^+-2}}=\lambda^{\frac{p^+}{p^+-2}}\geq\alpha^{p^+},\ \forall\ x\in\overline{\Omega}.
	\end{equation}
	
	\noindent This shows that $\displaystyle\int_{\Omega}\dfrac{z_n^{p(x)}}{\lambda^{\frac{p(x)}{p(x)-2}}}\ dx\leq\int_{\Omega} \dfrac{z_n^{p(x)}}{\alpha^{p^+}}\ dx\leq 1, \ \forall\ n\geq n_0$, which means in particular that $\lambda_1:=\alpha^{p^+-2}\in Z_n,\ \forall\ n\geq n_0$. Therefore:
	
	\begin{equation}
		\Vert z_n^{p(x)-2}\Vert_{L^{\frac{p(x)}{p(x)-2}}(\Omega)}=\inf Z_n\leq \lambda_1, \ \forall\ n\geq n_0,
	\end{equation}
	
	\noindent so that $\left (\Vert z_n^{p(x)-2}\Vert_{L^{\frac{p(x)}{p(x)-2}}(\Omega)} \right )_{n\geq 1}$ is bounded, as needed. This completes the proof of the fact that $\partial\mathcal{C}$ is a continuous operator. Using now Theorem \ref{athmgatfre} we conclude that $\mathcal{C}\in C^1\big (W^{1,p(x)}(\Omega) \big )$ and:
	
	\begin{equation}\label{3eqc3}
		\langle\mathcal{C}'(V),\phi\rangle=\int_{\Omega} |V(x)|^{p(x)-2}V(x)\phi(x)\ dx,\ \forall\ V,\phi\in W^{1,p(x)}(\Omega).
	\end{equation}

	\noindent Since $\mathcal{G}$ is a bounded linear functional it is continuously Fr\'{e}chet differentiable on all $W^{1,p(x)}(\Omega)$ and: 
	
	\begin{equation}\label{3eqc4}\langle \mathcal{G}'(V),\phi\rangle=\mathcal{G}(\phi)=\displaystyle\int_{\Omega} g(x)\phi(x)\ dx\ \text{for all}\ V,\phi\in W^{1,p(x)}(\Omega).
	\end{equation}

	\noindent Combining relations \eqref{3eqc1},\eqref{3eqc2},\eqref{3eqc3} and \eqref{3eqc4} we can conclude that $\mathcal{J}_{\lambda,\varepsilon}\in C^1\big (W^{1,p(x)}(\Omega) \big )$ and formula \eqref{3eqc} holds.
	
	\bigskip
	
	\noindent $\blacktriangleright$ \textbf{Fact III: $\mathcal{J}_{\lambda,\varepsilon}$ is a strictly convex functional.}
	
	\noindent Indeed, for every $V_1,V_2\in W^{1,p(x)}(\Omega)$ and $\theta\in [0,1]$, using the convexity of $A(x,\cdot)$ (see Proposition \ref{4prop1} \textbf{(4)}), we have that:
	
	\begin{align*}
		&\ \ \ \ \ \theta \mathcal{J}_{\lambda,\varepsilon}(V_1)+(1-\theta)\mathcal{J}_{\lambda,\varepsilon}(V_2)-\mathcal{J}_{\lambda,\varepsilon}\big (\theta V_1+(1-\theta V_2\big )=\\
		&=\int_{\Omega} \theta A(x,\nabla V_1)+(1-\theta)A(x,\nabla V_2)-A(x,\theta\nabla V_1+(1-\theta)\nabla V_2)\ dx+\dfrac{\lambda}{2}\theta(1-\theta)\int_{\Omega} (V_1-V_2)^2\ dx \\
		&+\varepsilon\int_{\Omega}\dfrac{1}{p(x)}\left (\theta|V_1|^{p(x)}+(1-\theta)|V_2|^{p(x)}-|\theta V_1+(1-\theta)V_2|^{p(x)} \right )\ dx\\
		&\geq\dfrac{\lambda}{2}\theta(1-\theta)\int_{\Omega} (V_1-V_2)^2\ dx
	\end{align*}
	
	\noindent Here we have used the fact that for each $x\in\overline{\Omega}$ the function $[0,\infty)\ni s\mapsto s^{p(x)}$ is strictly convex, because its first derivative $[0,\infty)\ni s\mapsto p(x)s^{p(x)-1}$ is a strictly increasing function as $p(x)\geq p^{-}>1$ for all $x\in\overline{\Omega}$. Thus we may write for a.a. $x\in\Omega$ that:
	
	\begin{equation}
		\theta|V_1|^{p(x)}+(1-\theta)|V_2|^{p(x)}\geq \left (\theta|V_1|+(1-\theta)|V_2| \right )^{p(x)}\stackrel{p(x)>1}{\geq} |\theta V_1+(1-\theta)V_2 |^{p(x)}.
	\end{equation}
	
	\noindent This shows that $\mathcal{J}_{\lambda,\varepsilon}$ is strictly convex.
	
	\bigskip
	
	\noindent $\blacktriangleright$ \textbf{Fact IV: $\mathcal{J}_{\lambda,\varepsilon}$ is a weakly lower semicontinuous functional.}
	
	\noindent Let now some $V\in W^{1,p(x)}(\Omega)$ and a sequence $(V_n)_{n\geq 1}\subset W^{1,p(x)}(\Omega)$ such that $V_n\rightharpoonup V$ in $W^{1,p(x)}(\Omega)$. Since $\mathcal{J}_{\lambda,\varepsilon}\in C^1\big (W^{1,p(x)}(\Omega)\big)$ we have that $\mathcal{J}'_{\lambda,\varepsilon}(V)\in W^{1,p(x)}(\Omega)^*$. Therefore: $\lim\limits_{n\to\infty} \langle \mathcal{J}'_{\lambda,\varepsilon}(V),V_n-V\rangle=0$.
	
	\noindent Using now the convexity of $\mathcal{J}_{\lambda,\varepsilon}$ and the fact that $\mathcal{J}_{\lambda,\varepsilon}\in C^1\big (W^{1,p(x)}(\Omega)\big)$ we have from \cite[Proposition 42.6, page 247]{zeidler} that for each $n\geq 1$:
	
	\begin{equation*}
		\mathcal{J}_{\lambda,\varepsilon}(V_n)\geq \mathcal{J}_{\lambda,\varepsilon}(V)+\langle \mathcal{J}'_{\lambda,\varepsilon}(V),V_n-V\rangle.
	\end{equation*}
	
	\noindent So, $\liminf\limits_{n\to\infty} \mathcal{J}_{\lambda,\varepsilon}(V_n)\geq \mathcal{J}_{\lambda,\varepsilon}(V)+\liminf\limits_{n\to\infty}\langle\mathcal{J}'_{\lambda,\varepsilon}(V),V_n-V\rangle=\mathcal{J}_{\lambda,\varepsilon}(V)$. This shows that $\mathcal{J}_{\lambda,\varepsilon}$ is weakly lower semicontinuous.
	
	\bigskip
	
	\noindent $\blacktriangleright$ \textbf{Fact V: $\mathcal{J}_{\lambda,\varepsilon}$ is a coercive functional for any $\varepsilon>0$.}
	
	\noindent Let $(V_n)_{n\geq 1}\subset W^{1,p(x)}(\Omega)$ with $\Vert V_n\Vert_{W^{1,p(x)}(\Omega)}\longrightarrow\infty$. From Theorem \ref{athsob} \textbf{(4)} we have that $\varrho_{p(x)}(V_n)=\rho_{p(x)}(V_n)+\rho_{p(x)}(|\nabla V_n|)\longrightarrow\infty$. From \textbf{(H11)}, we get that:
	
	\begin{align}
		\mathcal{J}_{\lambda,\varepsilon}(V_n)&=\int_{\Omega} A(x,\nabla V_n)\ dx+\dfrac{\lambda}{2}\int_{\Omega} V_n^2\ dx+\varepsilon\int_{\Omega}\dfrac{|V_n|^{p(x)}}{p(x)}\ dx-\int_{\Omega} gV_n\ dx\nonumber\\
		&\geq \delta\int_{\Omega}|\nabla V_n|^{p(x)}\ dx-\tilde{\delta}+\dfrac{\varepsilon}{p^+}\int_{\Omega} |V_n|^{p(x)}\ dx+\int_{\Omega} \dfrac{\lambda}{2}V_n^2-\lambda |V_n|\ dx\nonumber\\
		&\geq \min\left \{\delta,\dfrac{\varepsilon}{p^+}\right \}\varrho_{p(x)}(V_n)-\tilde{\delta}-\dfrac{\lambda}{2}|\Omega|\stackrel{n\to\infty}{\longrightarrow}\infty.
	\end{align}

	\noindent Combining the fact that $W^{1,p(x)}(\Omega)$ is a reflexive Banach space, \textbf{Fact III} and \textbf{Fact IV}, we deduce based on Theorem \ref{athmstru}, that there is some $V_{\varepsilon}\in W^{1,p(x)}(\Omega)$ such that: $\mathcal{J}_{\lambda,\varepsilon}(V_{\varepsilon})=\displaystyle\inf_{U\in W^{1,p(x)}(\Omega)} \mathcal{J}_{\lambda,\varepsilon}(U)$.
	
	\noindent Since $\mathcal{J}_{\lambda,\varepsilon}\in C^1\big (W^{1,p(x)}(\Omega) \big )$ and $V_{\varepsilon}$ is a (global) minimum of $\mathcal{J}_{\lambda,\varepsilon}$, we deduce that $\mathcal{J}_{\lambda,\varepsilon}'(V_{\varepsilon})=0$, i.e. $V_{\varepsilon}$ is a weak solution of the problem \eqref{eqelambdaeps}.\footnote{See Corollary 2.5 at page 53 from \cite{coleman}.}
\end{proof}

\begin{remark} For $\varepsilon=0$ we shall write $\mathcal{J}_{\lambda,\varepsilon}=\mathcal{J}_{\lambda,0}:=\mathcal{J}_{\lambda}$.	
\end{remark}

\begin{lemma}\label{3lemgama} Let any $V\in W^{1,p(x)}(\Omega)$ and any $(V_n)_{n\geq 1}\subset W^{1,p(x)}(\Omega)$ such that $V_n\longrightarrow V$ in $W^{1,p(x)}(\Omega)$. Then for every sequence of real numbers $(\varepsilon_n)_{n\geq 1}$ convergent to $0$ we have that:
	
	\begin{equation}
		\lim\limits_{n\to\infty} \mathcal{J}_{\lambda,\varepsilon_n}(V_n)=\mathcal{J}_{\lambda}(V).
	\end{equation}
	
\end{lemma}

\begin{proof} We have for each $n\geq 1$ that:
	
	\begin{align*} \big |\mathcal{J}_{\lambda,\varepsilon_n}(V_n)-\mathcal{J}_{\lambda}(V)\big |&=\left |\mathcal{J}_{\lambda}(V_n)-\mathcal{J}_{\lambda}(V)+\varepsilon_n\int_{\Omega}\dfrac{|V_n|^{p(x)}}{p(x)}\ dx\right |\\
		&\leq \big |\mathcal{J}_{\lambda}(V_n)-\mathcal{J}_{\lambda}(V) \big |+|\varepsilon_n|\int_{\Omega}\dfrac{|V_n|^{p(x)}}{p(x)}\ dx\\
		&\leq \big |\mathcal{J}_{\lambda}(V_n)-\mathcal{J}_{\lambda}(V) \big |+\dfrac{|\varepsilon_n|}{p^{-}}\int_{\Omega}|V_n|^{p(x)}\ dx\\
		&=\big |\mathcal{J}_{\lambda}(V_n)-\mathcal{J}_{\lambda}(V) \big |+\dfrac{|\varepsilon_n|}{p^{-}}\rho_{p(x)}(V_n)\\
		&\longrightarrow 0\ \text{as}\ n\to\infty.
	\end{align*}
	
	\noindent Above we have used the continuity of $\mathcal{J}_{\lambda}$ on $W^{1,p(x)}(\Omega)$ and the fact that $\big (\rho_{p(x)}(V_n)\big )_{n\geq 1}$ is a bounded sequence. Let's explain this fact in more details: from Theorem \ref{athleb} \textbf{(11)} we have that:
	
	\begin{equation}
		0\leq \rho_{p(x)}(V_n)=\rho_{p(x)}(V_n-V+V)\leq 2^{p^+-1}\rho_{p(x)}(V_n-V)+2^{p^+-1}\rho_{p(x)}(V).
	\end{equation}
	
	\noindent Since $V_n\to V$ in $W^{1,p(x)}(\Omega)\hookrightarrow L^{p(x)}(\Omega)$ we also have that $V_n\to V$ in $L^{p(x)}(\Omega)$. Therefore, based on Theorem \ref{athleb} \textbf{(7)} we have that $\lim\limits_{n\to\infty}\rho_{p(x)}(V_n-V)=0$. Hence $\big (\rho_{p(x)}(V_n-V)\big )_{n\geq 1}$ is a bounded sequence (being convergent). This shows that $\big (\rho_{p(x)}(V_n)\big )_{n\geq 1}$ is also bounded and the proof is complete.

	\begin{remark} In particular, Lemma \ref{3lemgama} shows that the family of functionals $\big (\mathcal{J}_{\lambda,\varepsilon} \big )_{\varepsilon>0}$ is $\Gamma$-convergent to $\mathcal{J}_{\lambda}$.
	\end{remark}

\end{proof}

\begin{theorem}\label{3thmrot} Problem \eqref{eqelambda} admits a unique weak solution for any $g\in\mathcal{M}_{\lambda}(\Omega)$ that will be denoted $V\in\mathcal{U}$.
\end{theorem}

%

\begin{proof} \noindent From Theorem \ref{3thsoleps} we get that for each $\lambda>0$ the functional $\mathcal{J}_{\lambda}=\mathcal{J}_{\lambda,0}\in C^1\big (W^{1,p(x)}(\Omega) \big )$ and:
	
	\begin{equation}
		\langle \mathcal{J}'_{\lambda}(V),\phi\rangle=\int_{\Omega} \mathbf{a}(x,\nabla V(x))\cdot \nabla\phi(x)\ dx+\lambda\int_{\Omega} V(x)\phi(x)\ dx-\int_{\Omega} g(x)\phi(x)\ dx.
	\end{equation}
	
	\noindent Thus $V\in W^{1,p(x)}(\Omega)$ is a weak solution of the problem \eqref{eqelambda} if and only if $V$ is a critical point of $\mathcal{J}_{\lambda}$, i.e.: $\mathcal{J}_{\lambda}'(V)=0$.
	
	\bigskip

	\noindent\textbf{(Uniqueness)} This follows from the uniqueness part of the Theorem \ref{3thsoleps} for $\varepsilon=0$.
	
	\bigskip
	
	\noindent\textbf{(Existence)}  We now want to show that all the requirements of the Theorem \ref{athmstru} are satisfied for the functional $\mathcal{J}_{\lambda}:X:=W^{1,p(x)}(\Omega)\cap\mathcal{U}\to\mathbb{R}$.
	
	\begin{enumerate}
		\item[$\bullet$] $W^{1,p(x)}(\Omega)$ is a \textbf{reflexive Banach space}\footnote{See Remark \ref{aremsob} from the Appendix.} and $X$ is a \textbf{weakly closed} subset of $W^{1,p(x)}(\Omega)$. 
		
		\noindent It is easy to see that $X$ is a convex subset of $W^{1,p(x)}(\Omega)$. We now show that it is strongly closed. Indeed, let $V\in W^{1,p(x)}(\Omega)$ and $(V_n)_{n\geq 1}\subset X$ such that $V_n\longrightarrow V$ in $W^{1,p(x)}(\Omega)$. It follows that $V_n\longrightarrow V$ in $L^{p(x)}(\Omega)\hookrightarrow L^{p^-}(\Omega)$. Therefore there is a subsequence $[0,1]\ni V_{n_k}(x)\stackrel{k\to\infty}{\longrightarrow} V(x)$ for a.a. $x\in \Omega$.\footnote{See the Corollary from \cite[Page 234]{Jones}.} So $V(x)\in [0,1]$ for a.a. $x\in\Omega$ which means that $V\in W^{1,p(x)}(\Omega)\cap\mathcal{U}=X$.
		
		\noindent Using now \textit{Mazur's theorem}\footnote{See Theorem 3.3.18 and Corollary 3.3.22 in \cite{papa1}.} we conclude that $X$ is weakly closed.
		
		
		\item[$\bullet$] $\mathcal{J}_{\lambda}=\mathcal{J}_{\lambda,0}$ is \textbf{weakly lower semicontinuous}.
		
		\item[$\bullet$] $\mathcal{J}_{\lambda}:X\to\mathbb{R}$ is \textbf{coercive}, meaning that for $V\in X$: $\lim\limits_{\Vert V\Vert_{W^{1,p(x)}(\Omega)}\to \infty} \mathcal{J}_{\lambda}(V)=\infty$.
		
		\noindent Let $(V_n)_{n\geq 1}\subset X=W^{1,p(x)}(\Omega)\cap\mathcal{U}$ with $\Vert V_n\Vert_{W^{1,p(x)}(\Omega)}\longrightarrow\infty$. From Theorem \ref{athsob} \textbf{(4)} we have that $\varrho_{p(x)}(V_n)=\rho_{p(x)}(V_n)+\rho_{p(x)}(|\nabla V_n|)\longrightarrow\infty$. The key step here is to observe that for each $n\geq 1$:
		
		\begin{equation}
			\rho_{p(x)}(V_n)=\int_{\Omega} |V_n(x)|^{p(x)}\ dx\leq \int_{\Omega} 1\ dx=|\Omega|\ \Longrightarrow\ \rho_{p(x)}(|\nabla V_n|)=\int_{\Omega} |\nabla V_n|^{p(x)}\ dx\longrightarrow\infty.
		\end{equation}
		
		\noindent Therefore,

		\begin{align}
			\mathcal{J}_{\lambda}(V_n)&=\int_{\Omega} A(x,\nabla V_n)\ dx+\dfrac{\lambda}{2}\int_{\Omega} V_n^2\ dx-\int_{\Omega} gV_n\ dx\nonumber\\
			&\geq \delta\int_{\Omega}|\nabla V_n|^{p(x)}\ dx-\tilde{\delta}+\int_{\Omega} \dfrac{\lambda}{2}V_n^2-\lambda |V_n|\ dx\nonumber\\
			&\geq\delta\rho_{p(x)}(|\nabla V_n|)-\tilde{\delta}-\dfrac{\lambda}{2}|\Omega|\stackrel{n\to\infty}{\longrightarrow}\infty.
		\end{align}

	\end{enumerate}
	
	\noindent Having now all the requirements fulfilled we deduce that there is some $V\in X= W^{1,p(x)}(\Omega)\cap\mathcal{U}$ such that $\mathcal{J}_{\lambda}(V)=\displaystyle\min_{U\in X} \mathcal{J}_{\lambda}(U)$. In what follows, we prove a very important lemma:
	
	\begin{lemma}\label{lemlocglob} We have that: $m:=\displaystyle\inf_{U\in X} \mathcal{J}_{\lambda}(U)=\mathcal{J}_{\lambda}(V)=\inf_{U\in W^{1,p(x)}(\Omega)} \mathcal{J}_{\lambda}(U):=d$.
	\end{lemma}
	
	\begin{proof}

		\noindent Note that, since $\mathcal{J}_{\lambda}=\mathcal{J}_{\lambda,0}$ is bounded from below, we have that $m,d\in\mathbb{R}$. It is obvious that since $X\subset W^{1,p(x)}(\Omega)$ we have that $m=\displaystyle\inf_{U\in X} \mathcal{J}_{\lambda}(U)\geq\inf_{U\in W^{1,p(x)}(\Omega)} \mathcal{J}_{\lambda}(U)=d$. We need to show that the reverse inequality is also true: $d\geq m$.

		\noindent For each $\varepsilon>0$ we take $V_{\varepsilon}\in X=W^{1,p(x)}(\Omega)\cap\mathcal{U}$ the unique weak solution of the problem \eqref{eqelambdaeps}, i.e. the unique (global) minimum of the functional $\mathcal{J}_{\lambda,\varepsilon}$. First, since $V_{\varepsilon}\in X$ we have that:
		
		\begin{align*}
			\min_{U\in W^{1,p(x)}(\Omega)} \mathcal{J}_{\lambda,\varepsilon}(U)&=\mathcal{J}_{\lambda,\varepsilon}(V_\varepsilon)=\mathcal{J}_{\lambda}(V_{\varepsilon})+\varepsilon\int_{\Omega}\dfrac{|V_{\varepsilon}|^{p(x)}}{p(x)}\ dx\geq \mathcal{J}_{\lambda}(V_{\varepsilon})+0\\
			(V_{\varepsilon}\in X)\ \ \ \ \ \ \ 		&\geq \mathcal{J}_{\lambda}(V)=m.
		\end{align*}

		\noindent Secondly we can remark that:
		
		\begin{align*} \mathcal{J}_{\lambda,\varepsilon}(V_{\varepsilon})=\min_{U\in W^{1,p(x)}(\Omega)} \mathcal{J}_{\lambda,\varepsilon}(U)&\leq \mathcal{J}_{\lambda,\varepsilon}(V)=\mathcal{J}_{\lambda}(V)+\varepsilon\int_{\Omega}\dfrac{|V|^{p(x)}}{p(x)}\ dx\\
			(V\in \mathcal{U})\ \ \ \ \ \ \ \	&\leq \mathcal{J}_{\lambda}(V)+\dfrac{\varepsilon}{p^{-}}\int_{\Omega} 1^{p(x)}\ dx\\
			&=m+\dfrac{\varepsilon |\Omega|}{p^{-}}
		\end{align*}
		
		\noindent Thence:
		
		\begin{equation}
			m\leq \min_{U\in W^{1,p(x)}(\Omega)} \mathcal{J}_{\lambda,\varepsilon}(U)\leq m+\dfrac{\varepsilon |\Omega|}{p^{-}},\ \forall\ \varepsilon>0.
		\end{equation}
		
		\noindent This shows that:
		
		\begin{equation}\label{3eqmd}
			\lim\limits_{\varepsilon\to 0^+} \min_{U\in W^{1,p(x)}(\Omega)} \mathcal{J}_{\lambda,\varepsilon}(U)=m=\min_{U\in X} \mathcal{J}_{\lambda}(U),\ \text{or}\ \lim\limits_{\varepsilon\to 0^+}\mathcal{J}_{\lambda,\varepsilon}(V_{\varepsilon})=\mathcal{J}_{\lambda}(V).
		\end{equation}
		
		\noindent Since $d=\displaystyle\inf_{U\in W^{1,p(x)}(\Omega)} \mathcal{J}_{\lambda}(U)$ we can take a minimizing sequence $(U_n)_{n\geq 1}\subset W^{1,p(x)}(\Omega)$ with:
		
		\begin{equation}
			d\leq \mathcal{J}_{\lambda}(U_n)\leq d+\dfrac{1}{n},\ \forall\ n\geq 1.
		\end{equation}
		
		\noindent The key step comes here. For any $\varepsilon>0$ and any $n\geq 1$ we have that:
		
		\begin{align*}
			d\leq\mathcal{J}_{\lambda}(U_n)&\leq \mathcal{J}_{\lambda,\varepsilon}(U_n)=\mathcal{J}_{\lambda}(U_n)+\varepsilon\int_{\Omega} \dfrac{|U_n(x)|^{p(x)}}{p(x)}\ dx\\
			&\leq d+\dfrac{1}{n}+\dfrac{\varepsilon}{p^{-}}\int_{\Omega} |U_n(x)|^{p(x)}\ dx\\
			&=d+\dfrac{1}{n}+\dfrac{\varepsilon}{p^{-}}\rho_{p(x)}(U_n).
		\end{align*}
		
		\noindent Because $U_n\in W^{1,p(x)}(\Omega)\subset L^{p(x)}(\Omega)$ we have that $0\leq\rho_{p(x)}(U_n)<\infty$. So for $\varepsilon_n=\dfrac{p^{-}}{n\big (\rho_{p(x)}(U_n)+1\big )}\leq\dfrac{p^{-}}{n}$ we have that:
		
		\begin{equation}
			d\leq\mathcal{J}_{\lambda}(U_n)\leq \mathcal{J}_{\lambda,\varepsilon_n}(U_n)\leq d+\dfrac{2}{n},\ \forall\ n\geq 1.
		\end{equation}
		
		\noindent From the \textit{squeezing principle} we obtain that $\lim\limits_{n\to\infty} \mathcal{J}_{\lambda,\varepsilon_n}(U_n)=d$. Finally from \eqref{3eqmd}, since $\varepsilon_n\to 0^{+}$, and Lemma \ref{3lemgama} we get that:
		
		\begin{equation}
			d=\lim\limits_{n\to\infty} \mathcal{J}_{\lambda,\varepsilon_n}(U_n)\geq\lim\limits_{n\to\infty} \mathcal{J}_{\lambda,\varepsilon_n}(V_{\varepsilon_n})=\mathcal{J}_{\lambda}(V)=m.
		\end{equation}

	\end{proof}

	\noindent Since $\mathcal{J}_{\lambda}\in C^1\big (W^{1,p(x)}(\Omega) \big )$ and $V$ is a (global) minimum of $\mathcal{J}_{\lambda}$ as $\mathcal{J}_{\lambda}(V)=m=d=\displaystyle\inf_{U\in W^{1,p(x)}(\Omega)} \mathcal{J}_{\lambda}(U)$, we deduce that $\mathcal{J}_{\lambda}'(V)=0$, i.e. $V$ is a weak solution of the problem \eqref{eqelambda}.\footnote{See Corollary 2.5 at page 53 from \cite{coleman}.}
\end{proof}

\section{Steady-states}

\begin{definition} $\bullet$ We say that $U\in W^{1,p(x)}(\Omega)$ is a \textbf{weak solution} for \eqref{eqeg} if $1\geq U\geq 0$ a.e. on $\Omega$ and for any test function $\phi\in W^{1,p(x)}(\Omega)$ we have that:
	\begin{equation}
		\int_{\Omega} \mathbf{a}(x,\nabla U(x))\cdot \nabla\phi(x)\ dx=\int_{\Omega} f(x,U(x))\phi(x)\ dx.
	\end{equation}
	
	$\bullet$ We say that $U\in W^{1,p(x)}(\Omega)$ is a \textbf{weak subsolution} for \eqref{eqeg} if $1\geq U\geq 0$ a.e. on $\Omega$ and for any test function $\phi\in W^{1,p(x)}(\Omega)^+$ we have that:
	\begin{equation}
		\int_{\Omega} \mathbf{a}(x,\nabla U(x))\cdot \nabla\phi(x)\ dx\leq\int_{\Omega} f(x,U(x))\phi(x)\ dx.
	\end{equation}
	
	$\bullet$ We say that $U\in W^{1,p(x)}(\Omega)$ is a \textbf{weak supersolution} for \eqref{eqeg} if $1\geq U\geq 0$ a.e. on $\Omega$ and for any test function $\phi\in W^{1,p(x)}(\Omega)^+$ we have that:
	\begin{equation}
		\int_{\Omega} \mathbf{a}(x,\nabla U(x))\cdot \nabla\phi(x)\ dx\geq\int_{\Omega} f(x,U(x))\phi(x)\ dx.
	\end{equation}
\end{definition}

\begin{remark} It is easy to notice from \textnormal{\textbf{(H10)}} that $0$ is a \textbf{weak subsolution} of $\eqref{eqeg}$ and $1$ is a \textbf{weak supersolution} of $\eqref{eqeg}$. Indeed:
	
	\begin{equation*}
		\begin{cases}\displaystyle\int_{\Omega} \mathbf{a}(x,\nabla 0)\cdot \nabla\phi(x)\ dx=0\leq \int_{\Omega} f(x,0)\phi(x)\ dx,\ \forall\ \phi\in W^{1,p(x)}(\Omega)^+ \\[5mm]
			\displaystyle\int_{\Omega} \mathbf{a}(x,\nabla 1)\cdot \nabla\phi(x)\ dx=0\geq \int_{\Omega} f(x,1)\phi(x)\ dx,\ \forall\ \phi\in W^{1,p(x)}(\Omega)^+
		\end{cases}.
	\end{equation*}
\end{remark}

\begin{remark}\label{3remsource} From \textnormal{\textbf{(H11)}} we have that there is some $\lambda_0>0$ such that the function $[0,1]\ni s\mapsto f(x,s)+\lambda_0 s$ is a strictly increasing function for a.a. $x\in\Omega$. We introduce the function: $g:\overline{\Omega}\times\mathbb{R}\to\mathbb{R},\ g(x,s)=f(x,s)+\lambda_0 s$. 
	
	\noindent We have here a very important observation to do: For any $U\in\mathcal{U}$ we have that the function $\Omega\ni x\longmapsto g(x,U(x))$ is from $\mathcal{M}_{\lambda_0}(\Omega)$. Indeed, using the monotonicity of $g(x,\cdot)$ on $[0,1]$ and \textnormal{\textbf{(H9)}} we have that:
	
	\begin{align*}
		&0\leq U(x)\leq 1\ \text{for a.a.}\ x\in\Omega \ \Longrightarrow\\[3mm]
		&0\leq f(x,0)=g(x,0)\leq g(x,U(x))\leq g(x,1)=f(x,1)+\lambda_0\leq\lambda_0\ \text{for a.a.}\ x\in\Omega.
	\end{align*}
	
	\noindent This shows that indeed $g(\cdot, U)\in\mathcal{M}_{\lambda_0}(\Omega)$ for any $U\in\mathcal{U}$.
	
\end{remark}

\begin{definition} Remark \ref{3remsource} and Theorem \ref{3thmrot} allows us to define the following (nonlinear) operator:
	
	\begin{equation}\label{3eqdefk}
		\mathcal{K}:\mathcal{U}\to \mathcal{U},\ \mathcal{K}(U)=V,\ \text{where}\ \begin{cases}-\operatorname{div}\mathbf{a}(x,\nabla V)+\lambda_0 V=g(x,U(x)), & x\in\Omega\\[3mm] \dfrac{\partial V}{\partial\nu}=0, & x\in\partial\Omega \end{cases}.
	\end{equation}
	
\end{definition}

\begin{proposition} The (nonlinear) operator $\mathcal{K}:\mathcal{U}\to\mathcal{U}$ has the following properties:
	\begin{enumerate}
		\item[\textbf{(1)}] $\mathcal{K}(U)\in W^{1,p(x)}(\Omega)\cap\mathcal{U}$ for any $U\in\mathcal{U}$;
		
		\item[\textbf{(2)}] $\mathcal{K}$ is strictly monotone, i.e. for any $U_1,U_2\in\mathcal{U}$ with $U_1\leq U_2$ a.e. on $\Omega$ and $U_1\not\equiv U_2$ we have that $\mathcal{K}(U_1)\leq \mathcal{K}(U_2)$ a.e. on $\Omega$ and $\mathcal{K}(U_1)\not\equiv\mathcal{K}(U_2)$;
		
		\item[\textbf{(3)}] $\mathcal{K}$ is $1+\dfrac{\gamma}{\lambda_0}$--Lipschitz continuous with respect to the norm of $L^2(\Omega)$, i.e. for any $U_1,U_2\in\mathcal{U}$:
		
		\begin{equation}\label{3ineqlip}
			\Vert \mathcal{K}(U_2)-\mathcal{K}(U_1)\Vert_{L^2(\Omega)}\leq \left (1+\dfrac{\gamma}{\lambda_0} \right )\Vert U_2-U_1\Vert_{L^2(\Omega)}.
		\end{equation}
	\end{enumerate}
	
\end{proposition}

\begin{proof} \textbf{(1)} This is obvious since $\mathcal{K}(U)$ is a weak solution of \eqref{3eqdefk}. Hence $\mathcal{K}(U)\in W^{1,p(x)}(\Omega)\cap\mathcal{U}=X$.
	
	\bigskip
	
	\noindent\textbf{(2)} Denote $V_1=\mathcal{K}(U_1)$ and $V_2=\mathcal{K}(U_2)$. Therefore:
	
	\begin{equation}
		\begin{cases}-\operatorname{div}\mathbf{a}(x,\nabla V_2)+\lambda_0 V_2=g(x,U_2(x))\geq g(x,U_1(x))=-\operatorname{div}\mathbf{a}(x,\nabla V_1)+\lambda_0 V_1, & x\in\Omega\\[3mm] \dfrac{\partial V_2}{\partial\nu}-\dfrac{\partial V_1}{\partial\nu}=0, & x\in\partial\Omega \end{cases}.
	\end{equation}
	
	\noindent From the \textit{weak comparison principle} we have that $V_2\geq V_1$ a.e. on $\Omega$ and if $V_2\equiv V_1$ then $g(x,U_1(x))=g(x,U_2(x))$ for a.a. $x\in\Omega$ which is the same as $U_1(x)=U_2(x)$ ($g(x,\cdot)$ is strictly increasing on $[0,1]$) for a.a. $x\in\Omega$. So $U_1\equiv U_2$. Thus $V_1\not\equiv V_2$ and we are done.	
	
	\noindent\textbf{(3)} Keeping the same notations as in the proof of \textbf{(2)} we have that:
	
	\begin{equation}
		\int_{\Omega} \big (\mathbf{a}(x,\nabla V_2)-\mathbf{a}(x,\nabla V_1) \big )\cdot\nabla\phi \ dx+\lambda_0\int_{\Omega} (V_2-V_1)\phi\ dx=\int_{\Omega} \big (g(x,U_2)-g(x,U_1)\big )\phi\ dx,\ \forall\ \phi\in W^{1,p(x)}(\Omega).
	\end{equation}
	
	\noindent Choosing $\phi=V_2-V_1\in W^{1,p(x)}(\Omega)$ we obtain that:
	
	\begin{align*}
		\lambda_0 \Vert V_2-V_1\Vert^2_{L^2(\Omega)}&\leq \int_{\Omega}\big (\mathbf{a}(x,\nabla V_2)-\mathbf{a}(x,\nabla V_1) \big )\cdot \big (\nabla V_2-\nabla V_1\big ) \ dx+\lambda_0\int_{\Omega} (V_2-V_1)^2\ dx\\
		&=\int_{\Omega} \big (g(x,U_2)-g(x,U_1))(V_2-V_1)\ dx\leq \int_{\Omega} \big |g(x,U_2)-g(x,U_1)|\cdot |V_2-V_1|\ dx\\
		\text{\textbf{(H10)}}\ \ \ \ \ \ \ 		&\leq (\lambda_0+\gamma)\int_{\Omega}|U_2-U_1|\cdot |V_2-V_1|\ dx \\
		\text{(Cauchy ineq.)}\ \ \ \ \ \ \ \ &\leq (\lambda_0+\gamma)\Vert U_2-U_1\Vert_{L^2(\Omega)}\Vert V_2-V_1\Vert_{L^2(\Omega)}
	\end{align*}
	
	\noindent Therefore: $\Vert V_2-V_1\Vert_{L^2(\Omega)}\leq\left (1+\dfrac{\gamma}{\lambda_0} \right )\Vert U_2-U_1\Vert_{L^2(\Omega)}$ and the proof is now complete.
\end{proof}

\begin{theorem} Problem \eqref{eqeg} admits at least two weak solutions (that can be equal). If they are equal, then there is no other solution of \eqref{eqeg}.
\end{theorem}

\begin{proof}\noindent\textbf{(Existence)} We introduce the following recurrence:
	
	\begin{equation}
		\begin{cases}U_{n+1}=\mathcal{K}(U_n),\ n\geq 0\\[3mm] U_0\in\mathcal{U}\end{cases}.
	\end{equation}
	
	\noindent $\bullet$ If we take $U_0\equiv 1\in\mathcal{U}$ then $U_1=\mathcal{K}(U_0)\in\mathcal{U}\ \Rightarrow\ U_1\leq 1=U_0$ a.e. on $\Omega$. Using the monotony of $\mathcal{K}$ we get immediately by induction that $U_2=\mathcal{K}(U_1)\leq \mathcal{K}(U_0)=U_1$, ... , $U_{n+1}=\mathcal{K}(U_n)\leq \mathcal{K}(U_{n-1})=U_n$ for any $n\geq 1$. Also, from $1\equiv U_0\geq 0$ a.e. on $\Omega$ we obtain that $U_1=\mathcal{K}(U_0)\geq \mathcal{K}(0)\geq 0$ because $\mathcal{K}(0)\in\mathcal{U}$. Repeating the same process $U_2=\mathcal{K}(U_1)\geq \mathcal{K}(0)\geq 0$ and in general we get that $U_n\geq 0$ for any $n\geq 1$. Hence:
	
	\begin{equation}
		1\equiv U_0\geq U_1\geq U_2\geq\dots\geq U_n\geq U_{n+1}\dots\geq 0
	\end{equation}
	
	\noindent So, for a.a. $x\in\Omega$, the sequence of real numbers $\big (U_n(x)\big )_{n\geq 1}\subset [0,1]$ is decreasing and bounded. We conclude that it is convergent to some limit denoted by $\overline{U}(x)\in [0,1]$. We can state that $U_n\to\overline{U}$ pointwise a.e. on $\Omega$. Now, since $U_n=|U_n|\leq 1\in L^2(\Omega)$ ($\Omega$ is bounded) for any $n\geq 1$ we obtain from \textit{Lebesgue dominated convergence theorem} that $U_n\to \overline{U}$ in $L^2(\Omega)$.

	\noindent Using \eqref{3ineqlip} we get that $\mathcal{K}(U_n)\to\mathcal{K}(\overline{U})$ in $L^2(\Omega)$. But $\mathcal{K}(U_n)=U_{n+1}\to \overline{U}$ in $L^2(\Omega)$. From the uniqueness of the limit we conclude that $\mathcal{K}(\overline{U})=\overline{U}$. So $\overline{U}\in W^{1,p(x)}(\Omega)\cap \mathcal{U}=X$ and $\overline{U}$ is a weak solution of the problem \eqref{eqeg}.
	
	\begin{remark} If $1\not\equiv\mathcal{K}(1)$ then $1\gneq U_1\gneq U_2\gneq\dots\gneq U_{n}\gneq U_{n+1}\gneq\dots\gneq \overline{U}$. This follows from the strict monotony of $\mathcal{K}$.
	\end{remark}
	
	\noindent $\bullet$ If we set $U_0\equiv 0$ then in the same manner as above we will get that:
	
	\begin{equation}
		1\geq\dots U_{n+1}\geq U_n\geq\dots\geq U_2\geq U_1\geq U_0\equiv 0.
	\end{equation}
	
	\noindent Thus for a.a. $x\in\Omega$ the sequence $\big (U_n(x) \big )_{n\geq 1}\subset [0,1]$ is increasing and bounded, hence convergent to some $\underline{U}(x)\in [0,1]$. So $U_n\to\underline{U}$ pointwise a.e. on $\Omega$. From $|U_n|=U_n\leq 1\in L^2(\Omega)$ we deduce from the \textit{Lebesgue dominated convergence theorem} that $U_n\to\underline{U}$ in $L^2(\Omega)$. 
	
	\noindent Now from the continuity of $\mathcal{K}$ with respect to the $L^2(\Omega)$-norm we conclude that $U_{n+1}=\mathcal{K}(U_n)\to\mathcal{K}(\underline{U})$ in $L^2(\Omega)$. But since $U_{n+1}\to \underline{U}$ in $L^2(\Omega)$ we deduce that $\underline{U}=\mathcal{K}(\underline{U})$. Thus $\underline{U}\in W^{1,p(x)}(\Omega)\cap \mathcal{U}=X$ is also a weak solution of the problem \eqref{eqeg}.
	
	\begin{remark} If $0\not\equiv\mathcal{K}(0)$ then $0\lneq U_1\lneq U_2\lneq\dots\lneq U_{n}\lneq U_{n+1}\lneq\dots\lneq \underline{U}$. 
	\end{remark}
	
	\begin{remark} If $U\in X$ is any weak solution of the problem \eqref{eqeg} then: $\underline{U}(x)\leq U(x)\leq\overline{U}(x)$ for a.a. $x\in\Omega$. Indeed, from $U\in\mathcal{U}$, we have that: $0\leq U\leq 1$. Applying $\mathcal{K}$ we obtain: $\mathcal{K}(0)\leq\mathcal{K}(U)=U\leq \mathcal{K}(1)$. Repeating this process $n\geq 1$ times we get:
		
		\begin{equation}
			\mathcal{K}^n(0)\leq U\leq\mathcal{K}^n(1).
		\end{equation}
		
		\noindent Making $n\to\infty$ we get that $\underline{U}\leq U\leq\overline{U}$ a.e. on $\Omega$.
		
		\noindent We will refer to $\overline{U}$ as the \textbf{maximal solution} of \eqref{eqeg} and to $\underline{U}$ as the \textbf{minimal solution} of the same problem. Of course, if $\underline{U}=\overline{U}$ then problem \eqref{eqeg} has exactly one solution.
	\end{remark}
	
	\begin{remark} If $f(x,0)\equiv 0$ then $0\equiv\mathcal{K}(0)$ and therefore $\underline{U}\equiv 0$. This will be called the \textbf{trivial solution} of the problem \eqref{eqeg}.
	\end{remark}

	\begin{remark} If $U,\tilde{U}\in X$ are two weak solutions of \eqref{eqeg} then the following relation holds:\footnote{Starting with this relation I have proved the uniqueness of the steady-state for the Neumann laplacian in the paper \cite{max1}.}
		
		\begin{align}
			\int_{\Omega} \big (\Psi(x,|\nabla\tilde{U}|)-\Psi(x,|\nabla U|)\big )\nabla\tilde{U}\cdot\nabla U\ dx&=\int_{\Omega} g(x,\tilde{U})U-g(x,U)\tilde{U}\ dx\nonumber \\
			&=\int_{\Omega} f(x,\tilde{U})U-f(x,U)\tilde{U}\ dx
		\end{align}
		
	\end{remark}

\end{proof}

\section{Conclusion} 

\noindent In a future research paper we shall discuss about the number of solutions of \eqref{eqeg}. Unlocking problem \eqref{eqelambda} will allow to start a theory of quasilinear parabolic equations with Neumann boundary conditions similar to the one developed in \cite{Giaco1}.

\section{Appendix}

\noindent In this appendix we gather some of the applied results in the paper. We begin with some basic properties of the variable exponent Lebesgue and Sobolev spaces.

\begin{definition} For any measurable exponent $p:\Omega\to [1,\infty)$ we define:
	
	\begin{equation*}
		L^{p(x)}(\Omega):=\left\{f:\Omega\to\mathbb{R}\ \text{measurable}\ \big |\  \int_{\Omega}\left |f(x) \right |^{p(x)}\ dx<\infty\right\}.
	\end{equation*}
	
\end{definition}

\begin{theorem}\label{athleb} Set $\rho_{p(x)}(u)=\displaystyle\int_{\Omega} |u(x)|^{p(x)}\ dx$ for $u\in L^{p(x)}(\Omega)$. Then:\footnote{See Theorem 1.3 from \cite{Fan2}.}
	
	\begin{enumerate}
		\item[\textnormal{\bf{(1)}}] $\Vert u\Vert_{L^{p(x)}(\Omega)}=a\ \Longleftrightarrow\ \rho_{p(x)}\left (\dfrac{u}{a} \right)=1$.
		\item[\textnormal{\bf{(2)}}] $\Vert u\Vert_{L^{p(x)}(\Omega)}=1\ \Longleftrightarrow\ \rho_{p(x)}\left (u \right)=1$.
		\item[\textnormal{\bf{(3)}}] $\Vert u\Vert_{L^{p(x)}(\Omega)}>1\ \Longleftrightarrow\ \rho_{p(x)}\left (u \right)>1$.
		\item[\textnormal{\bf{(4)}}] $\Vert u\Vert_{L^{p(x)}(\Omega)}<1\ \Longleftrightarrow\ \rho_{p(x)}\left (u \right)<1$.
		\item[\textnormal{\bf{(5)}}] $\Vert u\Vert_{L^{p(x)}(\Omega)}>1\ \Longrightarrow\ \Vert u\Vert_{L^{p(x)}(\Omega)}^{p_{-}}\leq \rho_{p(x)}(u)\leq \Vert u\Vert_{L^{p(x)}(\Omega)}^{p_{+}}$.
		\item[\textnormal{\bf{(6)}}] $\Vert u\Vert_{L^{p(x)}(\Omega)}<1\ \Longrightarrow\ \Vert u\Vert_{L^{p(x)}(\Omega)}^{p_{+}}\leq \rho_{p(x)}(u)\leq \Vert u\Vert_{L^{p(x)}(\Omega)}^{p_{-}}$.
		\item[\textnormal{\bf{(7)}}] For $(u_n)_{n\geq 1}\subset L^{p(x)}(\Omega)$ and $u\in L^{p(x)}(\Omega)$ we have: $\Vert u_n-u\Vert_{L^{p(x)}(\Omega)}\longrightarrow 0\ \Longleftrightarrow\ \rho_{p(x)}(u_n-u)\longrightarrow 0$.\footnote{The proof can be found in \cite[Proposition 2.56 and Corollary 2.58, page 44]{Cruz}.}
		\item[\textnormal{\bf{(8)}}] For $(u_n)_{n\geq 1}\subset L^{p(x)}(\Omega)$ we have: $\Vert u_n\Vert_{L^{p(x)}(\Omega)}\longrightarrow \infty\ \Longleftrightarrow\ \rho_{p(x)}(u_n)\longrightarrow \infty$.
		\item[\textnormal{\bf{(9)}}] For any $u\in L^{p(x)}(\Omega)$ we have that: $\lim\limits_{n\to\infty}\Vert u\chi_{D_n}\Vert_{L^{p(x)}(\Omega)}=0$, where $(D_n)_{n\geq 1}$ is a sequence of measurable sets included in $\Omega$ with the property that $\lim\limits_{n\to\infty} |D_n|=0$.\footnote{See Theorem 1.13 from \cite{Fan2}.}
		\item[\textnormal{\bf(10)}] Let $u\in L^{p(x)}(\Omega)$ and $v:\Omega\to\mathbb{R}$ a measurable function such that $|u|\geq |v|$ a.e. on $\Omega$. Then $v\in L^{p(x)}(\Omega)$ and $\Vert u\Vert_{L^{p(x)}(\Omega)}\geq \Vert v\Vert_{L^{p(x)}(\Omega)}$.
		\item[\textnormal{\bf(11)}] $\rho_{p(x)}(u+v)\leq 2^{p^+-1}\big (\rho_{p(x)}(u)+\rho_{p(x)}(v) \big )$.
		\item[\textnormal{\bf(12)}] If $\lambda\in [0,1]$ we have $\rho_{p(x)}(\lambda u)\leq \lambda \rho(u)$ for each $u\in L^{p(x)}(\Omega)$. 
		\item[\textnormal{\bf(13)}] If $\lambda\geq 1$ we have $\rho_{p(x)}(\lambda u)\geq \lambda \rho(u)$ for each $u\in L^{p(x)}(\Omega)$. 
		\item[\textnormal{\bf(14)}] If $u,v:\Omega\to\mathbb{R}$ are two measurable functions and $|u|\geq |v|$ a.e. on $\Omega$, then $\rho_{p(x)}(u)\geq \rho_{p(x)}(v)$ and the inequality is strict if $|u|\not\equiv|v|$.
		\item[\textnormal{\bf(15)}] For any $u\in L^{p(x)}(\Omega)$, the function $[1,\infty)\ni\lambda\longmapsto \rho_{p(x)}\left (\dfrac{u}{\lambda} \right )$ is continuous and decreasing. Furthermore: $\lim\limits_{\lambda\to\infty} \rho_{p(x)}\left (\dfrac{u}{\lambda}\right )=0$ for every $u\in L^{p(x)}(\Omega)$.\footnote{See \cite[Proposition 2.7, page 17]{Cruz}.}
		\item[\textnormal{\bf(16)}] $\rho_{p(x)}:L^{p(x)}(\Omega)\to\mathbb{R}$ is a convex and continuous functional.
		\item[\textnormal{\bf{(17)}}] If $\Omega\subseteq\Omega_1\cup\Omega_2$ then $\Vert u\Vert_{L^{p(x)}(\Omega)}\leq \Vert u\Vert_{L^{p(x)}(\Omega_1)}+\Vert u\Vert_{L^{p(x)}(\Omega_2)}$ for any $u\in L^{p(x)}(\Omega_1\cup\Omega_2)$.
		\item[\textnormal{\bf{(18)}}] If $D\subset\Omega$ are two measurable sets then $\Vert u\chi_{D}\Vert_{L^{p(x)}(\Omega)}=\Vert u\Vert_{L^{p(x)}(D)}$ for any $u\in L^{p(x)}(\Omega)$.
		\item[\textnormal{\bf{(19)}}] $\rho_{p(x)}:L^{p(x)}(\Omega)\to\mathbb{R}$ is a lower semicontinuous and a weakly lower semicontinuous functional.\footnote{See \cite[Theorem 2.1.17 and Remark 2.1.18]{Hasto}.}
		\item[\textnormal{\bf{(20)}}] If $\lambda\in [1,\infty)$ then $\lambda^{p^-}\rho_{p(x)}(u)\leq\rho_{p(x)}(\lambda u)\leq \lambda^{p^+}\rho_{p(x)}(u)$ for any measurable function $u:\Omega\to\mathbb{R}$.
		\item[\textnormal{\bf{(21)}}] If $\lambda\in (0,1)$ then $\lambda^{p^-}\rho_{p(x)}(u)\geq\rho_{p(x)}(\lambda u)\geq \lambda^{p^+}\rho_{p(x)}(u)$ for any measurable function $u:\Omega\to\mathbb{R}$.
		\item[\textnormal{\bf{(22)}}] If $\Vert u\Vert_{L^{p(x)}(\Omega)}>1$ then: $\rho_{p(x)}(u)^{\frac{1}{p^+}}\leq \Vert u\Vert_{L^{p(x)}(\Omega)}\leq\rho_{p(x)}(u)^{\frac{1}{p^-}}$.
		\item[\textnormal{\bf{(23)}}] If $0<\Vert u\Vert_{L^{p(x)}(\Omega)}<1$ then: $\rho_{p(x)}(u)^{\frac{1}{p^+}}\geq \Vert u\Vert_{L^{p(x)}(\Omega)}\geq\rho_{p(x)}(u)^{\frac{1}{p^-}}$.
		\item[\textnormal{\bf{(24)}}] If $p:\Omega\to [1,\infty)$ is a measurable bounded exponents and $q:\Omega\to (0,\infty)$ is also measurable such that $q(x)\leq p(x)$ a.e. on $\Omega$ then for any $u\in L^{p(x)}(\Omega)$ we have that $|u|^{q(x)}\in L^{\frac{p(x)}{q(x)}}(\Omega)$.
		\item[\textnormal{\bf{(25)}}] If $\Vert u\Vert_{L^{p(x)}(\Omega)}\leq 1$ then $\rho_{p(x)}(u)\leq \Vert u\Vert_{L^{p(x)}(\Omega)}$ and if $\Vert u\Vert_{L^{p(x)}(\Omega)}\geq 1$ then $\rho_{p(x)}(u)\geq \Vert u\Vert_{L^{p(x)}(\Omega)}$.\footnote{This is Corollary 2.22 from \cite[page 24]{Cruz}.}
		
	\end{enumerate}	
\end{theorem}

\begin{definition} For a measurable exponent $p:\Omega\to [1,\infty)$ let:
	\begin{align*}
		W^{1,p(x)}(\Omega)&=\left\{u\in L^{p(x)}(\Omega)\ |\ |\nabla u|\in L^{p(x)}(\Omega)\right\}\ \text{and:}\\
		\Vert u\Vert_{W^{1,p(x)}(\Omega)}&=\Vert u\Vert_{L^{p(x)}(\Omega)}+\Vert |\nabla u|\Vert_{L^{p(x)}(\Omega)},\ \forall\ u\in W^{1,p(x)}(\Omega).
	\end{align*}
	
	\noindent We also set $W^{1,p(x)}_0(\Omega)$ as the closure of $C^{\infty}_0(\Omega)$ in $\big (W^{1,p(x)}(\Omega),\Vert \cdot\Vert_{W^{1,p(x)}(\Omega)}\big )$. We can use in $W^{1,p(x)}_0(\Omega)$ the equivalent norm $W^{1,p(x)}_0(\Omega)\ni u\longmapsto\Vert |\nabla u|\Vert_{L^{p(x)}(\Omega)}$.
\end{definition}

\begin{theorem}\label{athsob} Set $\varrho_{p(x)}(u)=\displaystyle\int_{\Omega} |u(x)|^{p(x)}+|\nabla u(x)|^{p(x)}\ dx$ for $u\in W^{1,p(x)}(\Omega)$. Then:
	
	\begin{enumerate}
		\item[\textnormal{\bf{(1)}}] $\Vert u\Vert_{W^{1,p(x)}(\Omega)}>1\ \Longrightarrow\ \Vert u\Vert_{W^{1,p(x)}(\Omega)}^{p_{-}}\leq \varrho_{p(x)}(u)\leq \Vert u\Vert_{W^{1,p(x)}(\Omega)}^{p_{+}}$.
		\item[\textnormal{\bf{(2)}}] $\Vert u\Vert_{W^{1,p(x)}(\Omega)}<1\ \Longrightarrow\ \Vert u\Vert_{W^{1,p(x)}(\Omega)}^{p_{+}}\leq \varrho_{p(x)}(u)\leq \Vert u\Vert_{W^{1,p(x)}(\Omega)}^{p_{-}}$.
		\item[\textnormal{\bf{(3)}}] $\Vert u_n-u\Vert_{W^{1,p(x)}(\Omega)}\longrightarrow 0\ \Longleftrightarrow\ \varrho_{p(x)}(u_n-u)\longrightarrow 0$.
		\item[\textnormal{\bf{(4)}}] $\Vert u_n\Vert_{W^{1,p(x)}(\Omega)}\longrightarrow \infty\ \Longleftrightarrow\ \varrho_{p(x)}(u_n)\longrightarrow \infty$.
	\end{enumerate}
	
\end{theorem}

\begin{remark}\label{aremsob} For $1<p^{-}\leq p^{+}<\infty$ we have that $\big (W^{1,p(x)}(\Omega),\Vert \cdot\Vert_{W^{1,p(x)}(\Omega)}\big )$ and $W^{1,p(x)}_0(\Omega)$ are uniform convex\footnote{See \cite[Theorem 7]{Mendez}.} (hence reflexive) and separable\footnote{See \cite[Theorem 2.1]{Fan2} or \cite[Theorem 3.1]{kovacik}.} Banach spaces. Moreover, if $p$ is log-H\"{o}lder continuous and $\Omega$ is a Lipschitz bounded domain then $C^{\infty}(\Omega)$ is dense in $W^{1,p(x)}(\Omega)$ and $W^{1,p(x)}_0(\Omega)=W^{1,p(x)}(\Omega)\cap W^{1,1}_0(\Omega)$.\footnote{See \cite[Theorem 2.6]{Fan2}.}
\end{remark}

\begin{theorem}[\textbf{Embeddings}]\label{athmemb} Let $\Omega$ be a bounded Lipschitz domain and $p,q:\overline{\Omega}\to [1,\infty)$ be two continuous exponents with:\footnote{The proof can be seen in \cite[Theorem 2.3]{Fan2}. Also take a look at \cite[Theorem 2.2(c)]{Dinca2}.}
	\begin{equation*}
		p(x)<N\ \text{and}\ q(x)<p^*(x),\ \forall\ x\in\overline{\Omega}.
	\end{equation*}
	\noindent Then $W^{1,p(x)}(\Omega)\stackrel{\text{c}}{\hookrightarrow} L^{q(x)}(\Omega)$. This result is still true if we remove the condition $p^+<N$.\footnote{See the proof of Proposition 2.2 from \cite{Fan3}.}
\end{theorem}

\begin{theorem}Let $\Omega$ be a bounded Lipschitz domain and $p:\overline{\Omega}\to (1,\infty)$, $p\in\mathcal{P}^{\text{log}}(\Omega)$ be a log-H\"{o}lder continuous exponent. If $q:\Omega\to [1,\infty)$ is a bounded measurable exponent with:
	\begin{equation*}
		1\leq q(x)\leq p^*(x),\ \text{for a.a.}\ x\in\Omega,
	\end{equation*}
	\noindent then $W^{1,p(x)}(\Omega)\hookrightarrow L^{q(x)}(\Omega)$. In particular $W^{1,p(x)}(\Omega)\hookrightarrow L^{p^*(x)}(\Omega)$.\footnote{See Proposition 2.1. and Proposition 2.2. from \cite{Fan3}.}
\end{theorem}

\begin{theorem}[\textbf{Nemytsky Operators}]\label{athnem} Let $f:\Omega\times\mathbb{R}\to\mathbb{R}$ be a Carath\'{e}odory function and $p:\Omega\to [1,\infty)$ be a measurable and bounded exponent. For each function $u\in L^{p(x)}(\Omega)$ we consider $\mathcal{N}_f(u):\Omega\to\mathbb{R},\ \mathcal{N}_f(u)(x)=f(x,u(x)),\ x\in\Omega$. 
	
	\noindent $\bullet$ If there is a non-negative function $g\in L^{q(x)}(\Omega)$, where $q:\Omega\to [1,\infty)$ is a measurable bounded exponent, and a constant $c\geq 0$ such that:
	\begin{equation}\label{aeqnem}
		|f(x,s)|\leq g(x)+c|s|^{\frac{p(x)}{q(x)}},\ \forall\ s\in\mathbb{R}\ \text{and}\ \text{for a.a.}\ x\in\Omega
	\end{equation}
	\noindent then $\mathcal{N}_f:L^{p(x)}(\Omega)\to L^{q(x)}(\Omega)$ is a continuous and bounded operator.

	\noindent $\bullet$ Conversely if it happens that $\mathcal{N}_f:L^{p(x)}(\Omega)\to L^{q(x)}(\Omega)$, then $\mathcal{N}_f$ is continuous and bounded and moreover there is a non-negative function $g\in L^{q(x)}(\Omega)$ and some constant $c\geq 0$ such that \eqref{aeqnem} holds.\footnote{For a complete proof see \cite[Theorem 1.16]{Fan2} or \cite[Theorem 4.1 and Theorem 4.2]{kovacik}.}
\end{theorem}

\begin{remark}\label{arenem} Theorem \ref{athnem} generalizez to Carath\'{e}odory functions of the type $f:\Omega\times\mathbb{R}^N\to\mathbb{R}$. For each $u=(u_1,u_2,\dots,u_N)\in L^{p_1(x)}(\Omega)\times L^{p_2(x)}(\Omega)\times\dots\times L^{p_N(x)}(\Omega)$, where $p_1,p_2,\dots,p_N:\Omega\to [1,\infty)$ are measurable and bounded exponents, we take the \textit{Nemytsky operator} $\mathcal{N}_f(u):\Omega\to\mathbb{R},\ \mathcal{N}_f(u)(x)=f\big (x,u(x)\big )=f\big (x,u_1(x),u_2(x),\dots,u_N(x)\big )$.
	
	\noindent $\bullet$ If there is a non-negative function $g\in L^{q(x)}(\Omega)$, where $q:\Omega\to [1,\infty)$ is a measurable bounded exponent, and some constants $c_1,c_2,\dots,c_N\geq 0$ such that:
	\begin{equation}\label{aeqnem2}
		|f(x,s_1,s_2,\dots,s_N)|\leq g(x)+c_1|s_1|^{\frac{p_1(x)}{q(x)}}+\dots+c_N|s_N|^{\frac{p_N(x)}{q(x)}},\ \forall\ s_1,s_2,\dots,s_N\in\mathbb{R}\ \text{and}\ \text{for a.a.}\ x\in\Omega
	\end{equation}
	\noindent then $\mathcal{N}_f:L^{p_1(x)}(\Omega)\times L^{p_2(x)}(\Omega)\times\dots\times L^{p_N(x)}(\Omega)\to L^{q(x)}(\Omega)$ is a continuous and bounded operator.
	
	\noindent $\bullet$ Conversely if it happens that $\mathcal{N}_f:L^{p_1(x)}(\Omega)\times L^{p_2(x)}(\Omega)\times\dots\times L^{p_N(x)}(\Omega)\to L^{q(x)}(\Omega)$, then $\mathcal{N}_f$ is continuous and bounded and moreover there is a non-negative function $g\in L^{q(x)}(\Omega)$ and some constants $c_1,c_2,\dots,c_N\geq 0$ such that \eqref{aeqnem2} holds.
\end{remark}

\begin{theorem}\label{apthplus} Let $\Omega$ be a bounded domain, $p:\Omega\to [1,\infty)$ be a measurable exponent and $u\in W^{1,p(x)}(\Omega)$. Set $u^+=\max\{u,0\}$ and $u^{-}=-\min\{u,0\}$. Then $u^{+},u^{-}\in W^{1,p(x)}(\Omega)^+$ and moreover:\footnote{The same result holds for the space $W^{1,p(x)}_0(\Omega)$ with the stronger assumption that $p\in\mathcal{P}^{\text{log}}(\Omega)$. For a proof, see \cite[Lemma 3.3]{Rad1}.}
	
	\begin{equation}\label{appeqnab}
		\nabla u^{+}=\begin{cases}0,\ \text{a.e. on}\ \{x\in\Omega\ |\ u(x)\leq 0\}\\[3mm]\nabla u,\ \text{a.e. on}\ \{x\in\Omega\ |\ u(x)>0\} \end{cases},\ \nabla u^{-}=\begin{cases}-\nabla u, \ \text{a.e. on}\ \{x\in\Omega\ |\ u(x)<0\}\\[3mm] 0, \ \text{a.e. on}\ \{x\in\Omega\ |\ u(x)\geq 0\} \end{cases}
	\end{equation}
	
	\begin{equation}
		\text{and}\ \nabla |u|=\begin{cases}-\nabla u, \ \text{a.e. on}\ \{x\in\Omega\ |\ u(x)<0\}\\[3mm] 0,\ \text{a.e. on}\ \{x\in\Omega\ |\ u(x)= 0\}\\[3mm]\nabla u,\ \text{a.e. on}\ \{x\in\Omega\ |\ u(x)>0\} \end{cases}.
	\end{equation}
	
	\noindent As a consequence: $u=u^+-u^-,\ |u|=u^++u^-,\ \nabla u=\nabla u^+-\nabla u^{-}$, $u^+u^-=0=\nabla u^+\cdot\nabla u^-=0$ and $\big | \nabla |u|\big |=|\nabla u|$ a.e. on $\Omega$.
	
	\noindent Another important fact is that:
	
	\begin{equation}\label{aineqplusminus}
		\Vert u^+\Vert_{W^{1,p(x)}(\Omega)},\Vert u^-\Vert_{W^{1,p(x)}(\Omega)}\leq \Vert |u|\Vert_{W^{1,p(x)}(\Omega)}=\Vert u\Vert_{W^{1,p(x)}(\Omega)}.
	\end{equation}
\end{theorem}

\begin{theorem}\label{athmstru} Let $\big (X,\Vert\cdot\Vert_X\big )$ be a real reflexive Banach space and $M\subset X$ be a weakly closed subset of $X$. Suppose that the functional $\mathcal{J}:M\to\mathbb{R}\cup\{+\infty\}$ is:
	
	\begin{enumerate}
		\item[$\bullet$] \textbf{coercive on $M$}, i.e.: $\forall\ r>0,\ \exists\ R>0$ such that $\forall\ u\in M$ with $\Vert u\Vert_X\geq R$ it follows that $\mathcal{J}(u)\geq r$;
		
		\item[$\bullet$] \textbf{sequentially weakly lower semicontinuous on $M$}, i.e.: for any $u\in M$ and any sequence $(u_n)_{n\geq 1}\subset M$ such that $u_n\rightharpoonup u$ in $X$ there holds:
		
		\begin{equation*}
			\liminf\limits_{n\to\infty}\mathcal{J}(u_n)\geq \mathcal{J}(u).
		\end{equation*}
	\end{enumerate}
	
	\noindent Then $\mathcal{J}$ is bounded from below on $M$, meaning that $\displaystyle\inf_{u\in M} \mathcal{J}(u)>-\infty$, and $\mathcal{J}$ attains its infimum on $M$, i.e. $\exists\ u\in M$ such that $\mathcal{J}(u)=\displaystyle\inf_{v\in M} \mathcal{J}(v)$.\footnote{This is Theorem 1.2 that can be found at page 4 in \cite{struwe}.}
\end{theorem}

\begin{theorem}\label{athmgatfre} Let $\big (X,\Vert\cdot\Vert_X\big )$ and $\big (Y,\Vert\cdot\Vert_Y\big )$ be two real normed linear spaces. If the function $f:X\to Y$ is \textbf{G\^ateaux differentiable} in an open neighborhood $U\subset X$ of a point $x_0\in X$ and the mapping $\partial f: U\to\mathcal{L}(X,Y)$ is continuous then $f$ is \textbf{Fr\'{e}chet differentiable} at $x_0$ and $f'(x_0)=\partial f(x_0)$.\footnote{This result is taken from \cite[Proposition 3.2.15, page 133]{drabek}. It can also be found as \cite[Theorem 1.1.3, page 3]{chang}.}
\end{theorem}

\begin{lemma}\label{alemhos} Let $a,b\in\mathbb{R}$ and $p>1$. Then $\lim\limits_{\theta\to 0}\dfrac{|a+b\theta|^p-|a|^p}{p\theta}=|a|^{p-2}ab$.
\end{lemma}

\begin{proof} If $a>0$ then for $|\theta|$ sufficiently small we also have that $a+b\theta>0$. Therefore, using \textit{L'H\^{o}spital Rule}, we get that:
	
	\begin{equation}
		\lim\limits_{\theta\to 0} \dfrac{|a+b\theta|^p-|a|^p}{p\theta}=\lim\limits_{\theta\to 0}\dfrac{(a+b\theta)^p-a^p}{p\theta}=\lim\limits_{\theta\to 0} (a+b\theta)^{p-1}b\stackrel{p>1}{=}a^{p-1}b=|a|^{p-2}ab.
	\end{equation}
	
	\noindent If $a<0$ then for $|\theta|$ sufficiently small we also have that $a+b\theta<0$. So:
	
	\begin{equation}
		\lim\limits_{\theta\to 0} \dfrac{|a+b\theta|^p-|a|^p}{p\theta}=\lim\limits_{\theta\to 0}\dfrac{(-a-b\theta)^p-(-a)^p}{p\theta}=\lim\limits_{\theta\to 0} (-a-b\theta)^{p-1}(-b)\stackrel{p>1}{=}(-a)^{p-1}(-b)=|a|^{p-2}ab.
	\end{equation}
	
	\noindent For $a=0$:
	
	\begin{equation}
		\lim\limits_{\theta\to 0} \dfrac{|a+b\theta|^p-|a|^p}{p\theta}=\lim\limits_{\theta\to 0}\dfrac{|b|^p\cdot |\theta|^p}{p\theta}=\dfrac{|b|^p}{p}\lim\limits_{\theta\to 0}|\theta|^{p-1}\operatorname{sgn}(\theta)\stackrel{p>1}{=}0=|a|^{p-2}ab.
	\end{equation}
\end{proof}

\begin{lemma}\label{alplap} For every $\mathbf{u},\mathbf{v}\in\mathbb{R}^N$ the following inequalities hold for some constants $c_1,c_2>0$:\footnote{See Lemma 3.3 and Lemma 3.4 from \cite{bystrom} where you can find the best constants $c_1,c_2$ and the cases when equality holds.}
	
	\begin{equation} \big | |\mathbf{u}|^{p-2}\mathbf{u}-|\mathbf{v}|^{p-2}\mathbf{v}\big |\leq c_1 |\mathbf{u}-\mathbf{v}|^{p-1},\ \text{if}\ p\in (1,2),
	\end{equation}
	
	\noindent and:
	
	\begin{equation}
		\big | |\mathbf{u}|^{p-2}\mathbf{u}-|\mathbf{v}|^{p-2}\mathbf{v}\big |\leq c_2\big (|\mathbf{u}|+|\mathbf{v}| \big )^{p-2}\cdot |\mathbf{u}-\mathbf{v}|,\ \text{if}\ p\in (2,\infty).
	\end{equation}
	
\end{lemma}

	\bibliographystyle{apalike}
	\bibliography{elliptic}
\end{document}